\documentclass{article}
\usepackage[utf8]{inputenc}

\usepackage{graphicx,cite}
\usepackage{amsmath,amsfonts,latexsym,url}
\usepackage[version=4]{mhchem}
\usepackage{siunitx}
\usepackage{longtable,tabularx}
\usepackage{color,soul} % for \hl command (highlighting text)
\usepackage{subfig} % for multiple plots in one figure

\newcommand{\m}[1]{{\bf{#1}}}
\newcommand{\g}[1]{\boldsymbol #1}

\newcommand{\tr}{^{\sf T}}
\newcommand{\C}[1]{{\cal {#1}}}

\newcommand{\ddt}[2]{^{^ #1} \hspace{-1mm} \frac{d}{dt} \hspace{-0.5mm} \left({#2}\right)}
\newcommand{\ddtm}[2]{\hspace{-1mm} ^{^{^{^ #1}}} \hspace{-1.5mm} \frac{d}{dt} \left( {#2} \right)}
\newcommand{\om}[2]{{^\C{#1} \hspace{-0.1mm} \g{\omega} \hspace{-0.1mm} ^\C{#2}}}
\newcommand{\al}[2]{{^\C{#1} \hspace{-0.1mm} \g{\alpha} \hspace{-0.1mm} ^\C{#2}}}
\newcommand{\sk}{^{\textrm{x}}}
\newcommand{\br}[2]{{{\left\lbrace {#2} \right\rbrace}_{#1}}}
\newcommand{\brsk}[2]{{{\left\lbrace {#2} \right\rbrace}_{#1} \sk}}
\newcommand{\brtr}[2]{{{\left\lbrace {#2} \right\rbrace}_{#1} \tr}}
\newcommand{\w}[1]{\omega_{#1}}
\newcommand{\e}[1]{\epsilon_{#1}}
\newcommand{\edot}[1]{\dot{\epsilon}_{#1}}
\newcommand{\h}[1]{{\m{#1}}}

\newcommand{\lsup}[2]{{\vphantom{#2}}^{#1}{{#2}}}

\newcommand{\CBA}{\m{C}_{\C{B} \C{A}}}
\newcommand{\CAB}{\m{C}_{\C{A} \C{B}}}

\newcommand{\CAE}{\m{C}_{\C{A} \C{E}}}
\newcommand{\CEA}{\m{C}_{\C{E} \C{A}}}

\newcommand{\CBE}{\m{C}_{\C{B} \C{E}}}
\newcommand{\CEB}{\m{C}_{\C{E} \C{B}}}

\title{\bf Nonsingular Euler Parameterizations for \\ Motion of a Point Mass in Atmospheric Flight}

\author{Alexander T.~Miller\footnote{Ph.D.~Student, NDSEG Fellow, Department of Mechanical and Aerospace Engineering.  E-mail:  alexandertmiller@ufl.edu.} \\ Anil~V.~Rao\footnote{Professor, Erich Farber Faculty Fellow, and University Term Professor, Department of Mechanical and Aerospace Engineering.  E-mail:  anilvrao@ufl.edu.  Associate Fellow AIAA.  Corresponding Author.} \vspace*{12pt} \\ {\em University of  Florida} \\ {\em Gainesville, FL, 32611}}

\begin{document}

\date{}
\maketitle{}
% \doublespacing
\renewcommand{\baselinestretch}{1}\normalsize\normalfont 

%-------------------------
 \begin{abstract}
   Three parameterizations are developed for modeling translational motion of a point mass in atmosphere flight over a central rotating body.  Unlike well-known parameterizations such as spherical coordinate parameterizations, where position and velocity are parameterized using a magnitude an an Euler angle rotation sequence, the method presented in this research employs Euler parameters.  Consequently, singularities and trigonometric functions are eliminated from the differential equations of motion.  As a result, the new parameterizations presented in this paper offer computational advantages over standard parameterizations that employ Euler angle sequences.  Finally, an example is studied where an atmospheric vehicle moves while in vertical flight, demonstrating the nonsingular nature of the formulations developed in this paper.
 \end{abstract}

% ----------------------------------------------------------------
\section*{Nomenclature}
\renewcommand{\baselinestretch}{1}\normalsize\normalfont
\begin{tabular}{lcl}
Symbol           &$=$& Description \\
$\m{a}$             &$=$& acceleration vector \\
$\m{C}$            &$=$& direction cosine matrix \\
$e_r$                 &$=$& position error \\
ECEF           &$=$& Earth-centered, Earth-fixed \\
ECI              &$=$& Earth-centered inertial \\
$\C{F}$             &$=$& reference frame \\
$\m{F}$             &$=$& force vector \\
$\m{h}$             &$=$& angular momentum vector \\
$h$                    &$=$& angular momentum \\
$\m{I}$              &$=$& identity matrix \\
LEO             &$=$& low Earth orbit \\
LVLH           &$=$& local vertical, local horizontal \\
$m$                   &$=$& mass \\
$\C{O}$             &$=$& reference point \\
$\C{P}$             &$=$& particle or center of mass \\
$\m{P}$             &$=$& general matrix \\
$\m{p}$             &$=$& general vector \\
$p$                    &$=$& general scalar \\
$\h{q}$              &$=$& unit vector along axis of rotation \\
$R_{e}$      &$=$& Earth radius \\
$\m{r}$              &$=$& position vector \\
$r$                     &$=$& radius \\
$\dot{r}$           &$=$& rate of change of radius \\
$T$                    &$=$& thrust magnitude \\
$t$                     &$=$& time on time interval $t \in [t_0,t_f]$ \\
$t_0$                 &$=$& initial time \\
$t_f$                  &$=$& terminal time \\
$\m{v}$             &$=$& velocity vector \\
$v$                    &$=$& speed \\
$\dot{v}$           &$=$& rate of change of speed \\
$\g{\alpha}$      &$=$& angular acceleration \\
$\m{\epsilon}$  &$=$& unit quaternion vector part \\
$\dot{\m{\epsilon}}$ &$=$& rate of change of $\m{\epsilon}$ \\
$\eta$                &$=$& unit quaternion scalar part \\
$\dot{\eta}$       &$=$& rate of change of $\eta$ \\
$\mu_{e}$  &$=$& Earth gravitational parameter \\
$\phi$                &$=$& general angle of rotation \\
$\g{\omega}$    &$=$& angular velocity \\
$\w{e}$      &$=$& Earth rotation rate
\end{tabular}

\renewcommand{\baselinestretch}{2}\normalsize\normalfont

%--------------------------
\section{Introduction}

In recent years, renewed attention has been given toward developing coordinate systems to describe the translational motion of a point in three-dimensional Euclidean space.  A great deal of this research has emanated from the orbital mechanics community, focusing on the use of quaternions for regularizing the equations of motion.  While traditionally unit quaternions (also known as Euler parameters) are used to describe the orientation of a rigid body \cite{Shuster1993, Kane1983, Wie1985, Wie2008, Hughes2, Junkins2012, Wertz2012}, more recent work has focused on extending quaternions to describe translational motion (see Refs.~\cite{KS1965, Chelnokov1981, Vivarelli1983, Deprit1994, Vrbik1994, Vrbik1995, Waldvogel2006, Waldvogel2008, Saha2009, Broucke1971, Deprit1975, Gurfil2005, Pelaez2007, Bau2015, Bau2020, Chelnokov2001, Chelnokov2003, Chelnokov2013, Chelnokov2014, Chelnokov2019, Libraro2014, Roa2017, Andreis2004} and the references therein).  In particular, the seminal work of Ref.~\cite{KS1965} used spinors to model three-dimensional motion in four dimensions.  Although the work of Ref.~\cite{KS1965} recognized the similarity between their KS transformation and quaternion multiplication, it was not until the works of Refs.~\cite{Chelnokov1981, Vivarelli1983} that the connection was made between the KS transformation and quaternions.  The connection between the KS transformation and quaternions has since been revisited and expanded upon in the works of Refs.~\cite{Deprit1994, Vrbik1994, Vrbik1995, Waldvogel2006, Waldvogel2008, Saha2009}.  In addition to the KS transformation, quaternions have been introduced to other translational motion descriptions.  The first three angles in the set of classical orbital elements, longitude of the ascending node, orbital inclination, and argument of periapsis, form a $3-1-3$ Euler angle sequence which is equivalently expressed by Euler parameters to eliminate the singularity associated with an equatorial orbit \cite{Broucke1971, Deprit1975, Gurfil2005}.  Finally, Refs.~\cite{Pelaez2007, Bau2015, Bau2020, Chelnokov2001, Chelnokov2003, Chelnokov2013, Chelnokov2014, Chelnokov2019,Libraro2014,Roa2017,Andreis2004} have further developed parameterizations using quaternions for orbital motion.

The aforementioned research has focused on orbital mechanics.  As such, these works focus primarily (or solely) on position--dependent forces, namely gravity.  The focus on position--dependent forces generally leads to the development of parameterizations that employ a single quaternion for quantifying position (for example, Refs.~\cite{Chelnokov2001, Libraro2014}).  Atmospheric flight mechanics, on the other hand, requires both position-- and velocity--dependent forces be taken into account.  When aerodynamic forces such as lift and drag are included, quantities such as the relative velocity, speed, and a convenient set of coordinates for quantifying the aerodynamic forces become important.  Thus, atmospheric flight mechanics problems are better suited to parameterizations that employ two quaternions for parameterizing both position and relative velocity.

In this paper, three parameterizations are developed for modeling the motion of a point mass over a central rotating body in atmospheric flight using two sets of Euler parameters (unit quaternions).  Similar to standard parameterizations such as spherical coordinates, the parameterizations developed in this paper decouple position and velocity into separate variables for magnitude and orientation.  Different from standard parameterizations, however, where the orientation variables correspond to angles of rotation in a sequence of four principal rotations (for instance, longitude, latitude, azimuth, and flight path angle), in this paper a sequence of two sets of Euler parameters (termed here as $\CAE$ and $\CBA$ Euler parameters) are developed.  The $\CAE$ Euler parameters define the orientation of a position frame (essentially a non-traditional local vertical, local horizontal frame) relative to the observation frame of the central, rotating body.  Likewise, the $\CBA$ Euler parameters define the orientation of a velocity frame (essentially a non-traditional wind frame) relative to the position frame.  The use of Euler parameters as orientation variables allows the equations of motion to be written in a form where singularities are eliminated at a pole or while the vehicle is in vertical flight.

To the best knowledge of the authors, no previous work has employed a sequence of two sets of Euler parameters in a point mass 3DOF model.  In addition, the intricacies involved in incorporating key atmospheric flight variables, namely the angle of attack and bank angle, into a quaternion representation of three degree of freedom motion, has not been explored previously.  Of all the references previously mentioned, only Ref.~\cite{Andreis2004} employs a quaternion for velocity (in their example).  The derivation in Ref.~\cite{Andreis2004}, however, assumes a single, full quaternion in the 3DOF model, not a sequence of two unit quaternions.  Another difference is that Ref.~\cite{Andreis2004} does not explore the ambiguity in their parameterization introduced by one component of the angular velocity being arbitrary, whereas this paper explicitly defines the source of the ambiguity via the language of reference frames.  Moreover, the present work explores three unique parameterizations, each with distinct characteristics, by enforcing different constraints on the arbitrary angular velocity term.  Finally, it is noted that the choice of frames in the example of Ref.~\cite{Andreis2004} leads to a singularity during vertical flight (similar to the $rvh$-Euler parameterization of the present research).  In contrast, the $rv$- and $rvL$-Euler parameterizations derived in this paper are nonsingular in vertical flight.

%Previously, Ref.~\cite{Andreis2004} has also parameterized position and velocity using quaternions, thus, combining magnitude and orientation variables into a single full quaternion.  It is noted, however, that combining  magnitude and orientation in the manner described in Ref.~\cite{Andreis2004} potentially leads to scaling issues.  Equally importantly, the parameterization developed in Ref.~\cite{Andreis2004} is undefined in vertical flight.  Finally, similar to the aforementioned research, Ref.~\cite{Andreis2004} focuses solely on orbital motion.  Thus, to the best knowledge of the authors of this paper, the intricacies involved in incorporating key atmospheric flight variables, namely the angle of attack and bank angle, into a quaternion representation of three degree of freedom motion, has not been explored previously.

The remainder of the paper is organized as follows.  Section~\ref{sect:NC} details the notation and conventions which will be employed throughout the paper.  A synopsis of the relevant math is provided in Section~\ref{sect:Review}.  The derivation of the equations of motion are contained in Section~\ref{sect:Derivation}.  After the equations of motion have been derived, Section~\ref{sect:Example} demonstrates their ability to model vertical flight without singularities.  Finally, key aspects of the three parameterizations are discussed in Section~\ref{sect:Discussion} and Section~\ref{sect:Conclusion} contains conclusions of the research.

% ----------------------------------------------------------------
\section{Notation and Conventions \label{sect:NC}}
The following notation and conventions will be employed throughout the paper.  Scalars will be represented by lowercase symbols (for example, $p \in \mathbb{R}$).  Next, vectors lie in three-dimensional Euclidean space and will be denoted by lowercase bold symbols (for example, $\m{p} \in \mathbb{E}^3$).  Moreover, the usual notation ``$\cdot$" and ``$\times$" specify scalar and vector products.  

Next, reference frames will be represented by uppercase calligraphic letters (for example, $\C{A}$).  Each reference frame will be designated a single coordinate system fixed in the frame.  Moreover, each coordinate system will be comprised of an origin and a set of three, right-handed, orthonormal basis vectors.  Table~\ref{tab:frames} summarizes the notation for the various reference frames, origins, and basis vectors which will be used.  The precise definitions of the basis vectors will be laid out in Section~\ref{sect:Derivation}.

\begin{table}[h]
  \centering
  \caption{Reference Frames.\label{tab:frames}}
  \renewcommand{\baselinestretch}{1}\normalsize\normalfont
 \begin{tabular}{|c|c|c|c|}\hline
    \textbf{Frame} & \textbf{Notation} & \textbf{Basis} & \textbf{Origin} \\\hline
    Inertial             & $\C{N}$     & $\{ \h{n}_1, \h{n}_2, \h{n}_3 \}$ & $O$ \\
    Observation & $\C{E}$      & $\{ \h{e}_1, \h{e}_2, \h{e}_3 \}$   & $O$ \\
    Position           & $\C{A}$      & $\{ \h{a}_1, \h{a}_2, \h{a}_3 \}$   & $O$ \\
    Velocity           & $\C{B}$      & $\{ \h{b}_1, \h{b}_2, \h{b}_3 \}$   & $P$  \\\hline
  \end{tabular}
\end{table}

Next, consider the expressions of vectors in particular bases.  The notation $\br{\C{A}}{\m{p}}$ will be employed to denote the expression of a vector $\m{p}$ in the (unique) basis assigned to frame $\C{A}$.  Moreover, $\br{\C{A}}{\m{p}}$ is equivalent to the column matrix $[p_1~p_2~p_3]\tr \in \mathbb{R}^3$ if $\m{p} = p_1 \h{a}_1 + p_2 \h{a}_2 + p_3 \h{a}_3$.  Furthermore, all direction cosine matrices will be represented by an uppercase bold ``$\m{C}$" with a subscript denoting the relevant frames and indicating the direction of the transformation.  Note that the relevant basis vectors are implied by the frame as laid out in Table~\ref{tab:frames}.  Thus, the identities
\begin{equation}\label{eq:DCM-def}
\begin{array}{lclclcl}
\br{\C{B}}{\m{p}} &=& \CBA \br{\C{A}}{\m{p}} & \textrm{and} & \br{\C{A}}{\m{p}} &=& \CAB \br{\C{B}}{\m{p}},
\end{array}
\end{equation}
hold for any arbitrary vector $\m{p}$, choice of frames $\C{A}$ and $\C{B}$, and choice of right-handed, orthonormal bases $\{ \h{a}_1, \h{a}_2, \h{a}_3 \}$ for frame $\C{A}$ and $\{ \h{b}_1, \h{b}_2, \h{b}_3 \}$ for frame $\C{B}$.  It is noted that $\CAB = \CBA\tr = \CBA^{-1}$ as well.  Finally, the notations $\CAB(i,:)$, $\CAB(:,j)$, and $\CAB(i,j)$ will be used to denote, respectively, the $i^{th}$ row of $\CAB$, the $j^{th}$ column of $\CAB$, or element $(i,j)$ of $\CAB$.

Next, scalar and vector products of two vectors will be expressed in the basis of a desired frame as follows.  Let $\m{p}$ and $\m{q}$ be any two vectors and let $\C{A}$ be the desired frame.  Moreover, let $\br{\C{A}}{\m{p}} = [p_1~p_2~p_3]\tr$ and $\br{\C{A}}{\m{q}} = [q_1~q_2~q_3]\tr$ be the expressions of vectors $\m{p}$ and $\m{q}$ in the basis fixed in frame $\C{A}$. Then, the scalar product of $\m{p}$ with $\m{q}$ expressed in the basis fixed in frame $\C{A}$ will be denoted by $\br{\C{A}}{\m{p} \cdot \m{q}} = \brtr{\C{A}}{\m{p}} \br{\C{A}}{\m{q}}$.  Likewise, the vector product of $\m{p}$ with $\m{q}$ expressed in the basis fixed in frame $\C{A}$ will be represented by $\br{\C{A}}{\m{p} \times \m{q}} = \brsk{\C{A}}{\m{p}} \br{\C{A}}{\m{q}}$, where
\begin{equation}\label{eq:skew}
\brsk{\C{A}}{\m{p}} = \begin{bmatrix}
0      & -p_3 & p_2\\
p_3   & 0      & -p_1\\
-p_2 & p_1   & 0
\end{bmatrix},
\end{equation}
denotes the skew-symmetric operator $\{ \cdot \}\sk$ acting on $\br{\C{A}}{\m{p}}$.  Finally, the identities
\begin{equation}\label{eq:dot-cross-express}
\begin{array}{lclclcl}
\brtr{\C{B}}{\m{p}} \br{\C{B}}{\m{q}} &=& \brtr{\C{A}}{\m{p}} \br{\C{A}}{\m{q}} & \textrm{and} & \brsk{\C{B}}{\m{p}} \br{\C{B}}{\m{q}}&=& \CBA \brsk{\C{A}}{\m{p}} \br{\C{A}}{\m{q}},
\end{array}
\end{equation}
relate the expressions of the scalar and vector products of vectors $\m{p}$ and $\m{q}$ in the bases of any two frames $\C{A}$ and $\C{B}$.

Lastly, the rate of change of scalars and vectors are defined as follows. First, all rates of change of scalars are denoted by the overdot symbol.  For example, the notation $\dot{p} = dp/dt$ denotes the rate of change of the scalar $p$.  Next, all rates of change of vectors are denoted with a left superscript to indicate the reference frame in which the rate of change is observed.  For example, the notation $\lsup{\C{A}}{d\m{p}/dt}$ denotes the rate of change of the vector $\m{p}$ as viewed by an observer in frame $\C{A}$.  Moreover, the rate of change of the vector $\m{p}$ as viewed by an observer in frame $\C{A}$ is related to the rate of change of the vector $\m{p}$ as viewed by an observer in frame $\C{B}$ via the transport theorem $\lsup{\C{A}}{d\m{p}/dt} = \lsup{\C{B}}{d\m{p}/dt} + \om{A}{B} \times \m{p}$, where $\om{A}{B}$ is the angular velocity of frame $\C{B}$ as viewed by an observer in frame $\C{A}$.  Similarly, the angular acceleration of frame $\C{B}$ as viewed by an observer in frame $\C{A}$ is denoted by $\al{A}{B}$.  Finally, the velocity and acceleration as viewed by an observer in reference frame $\C{A}$ are denoted, respectively, as $\lsup{\C{A}}{\m{v}}$ and $\lsup{\C{A}}{\m{a}}$, where $\C{A}$ denotes the frame in which the rate of change is observed.  

% ----------------------------------------------------------------
\section{Mathematical Preliminaries \label{sect:Review}}
A brief review of Euler parameters is now presented that will be used later to develop the $rv$-Euler parameterization.  Consider a direction cosine matrix $\CBA$ where frames $\C{A}$ and $\C{B}$ (and their right-handed, orthonormal bases) are arbitrary.  It is known that a right-handed orthonormal basis $\{ \h{a}_1, \h{a}_2, \h{a}_3 \}$ fixed in frame $\C{A}$ can be rotated to be aligned with the basis $\{ \h{b}_1, \h{b}_2, \h{b}_3 \}$ fixed in frame $\C{B}$ by a single rotation about an axis $\h{q}$ by an angle $\phi$.  Suppose further that $\h{q}$ is a unit vector which, when expressed in the basis $\{ \h{b}_1, \h{b}_2, \h{b}_3 \}$, is given as $\br{\C{B}}{\h{q}} = [q_1, q_2, q_3]\tr$, noting that, for the specific case where $\m{q}$ is the axis of rotation, $\br{\C{B}}{\h{q}} = \br{\C{A}}{\h{q}}$.  Then, the axis-angle representation of $\CBA$ is written as
\begin{equation}\label{eq:C_BA-AxisAngle}
\CBA = \begin{bmatrix}
(1 - \cos \phi) q_1^2 + \cos \phi &
(1 - \cos \phi) q_1 q_2 + q_3 \sin \phi &
(1 - \cos \phi) q_1 q_3 - q_2 \sin \phi \\
(1 - \cos \phi) q_2 q_1 - q_3 \sin \phi &
(1 - \cos \phi) q_2^2 + \cos \phi &
(1 - \cos \phi) q_2 q_3 + q_1 \sin \phi \\
(1 - \cos \phi) q_3 q_1 + q_2 \sin \phi &
(1 - \cos \phi) q_3 q_2 - q_1 \sin \phi &
(1 - \cos \phi) q_3^2 + \cos \phi
\end{bmatrix}.
\end{equation}

The axis-angle parameters give rise to a physical interpretation of the Euler parameters which are defined from the axis-angle parameters as
\begin{equation}\label{eq:Quat}
\begin{array}{lclclclclclclcl}
\e{1} & = & q_1 \sin \phi/2 & , & \e{2} & = & q_2 \sin \phi/2 & , & \e{3} & = & q_3 \sin \phi/2 & , & \eta & = & \cos \phi/2.  
\end{array}
\end{equation}
In Eq.~\eqref{eq:Quat} it is observed that the Euler parameters $\{ \e{1}, \e{2}, \e{3} \}$ provide the same directional information as the axis-angle parameters $\{ q_1, q_2, q_3 \}$.  Note also that the Euler parameters satisfy $\e{1}^2 + \e{2}^2 + \e{3}^2 + \eta^2 = 1$
which is equivalent to a unit norm constraint on the quaternion defined by the vector part $[\e{1}, \e{2}, \e{3}]\tr$ and scalar part $\eta$.  Now, given the definitions for the Euler parameters in Eq.~\eqref{eq:Quat}, the description of $\CBA$ in Eq.~\eqref{eq:C_BA-AxisAngle} becomes
\begin{equation}\label{eq:C_BA-Quat}
\CBA = \begin{bmatrix}
1 - 2(\e{2}^2 + \e{3}^2) &
2(\e{1}\e{2} + \e{3} \eta) &
2(\e{1}\e{3} - \e{2} \eta) \\
2(\e{2}\e{1} - \e{3} \eta) &
1 - 2(\e{3}^2 + \e{1}^2) &
2(\e{2}\e{3} + \e{1} \eta) \\
2(\e{3}\e{1} + \e{2} \eta) &
2(\e{3}\e{2} - \e{1} \eta) &
1 - 2(\e{1}^2 + \e{2}^2)
\end{bmatrix}.
\end{equation}
In addition, the rates of change of the Euler parameters are related to the angular velocity $\om{A}{B}$ as
\begin{equation}\label{eq:QuatRate}
\begin{array}{lclclcl}
\dot{\e{}}_1 &=& \phantom{-} \frac{1}{2} \left( \phantom{-} \eta \w{1} - \e{3} \w{2} + \e{2} \w{3} \right) & , & 
\dot{\e{}}_2 &=& \phantom{-} \frac{1}{2} \left( \e{3} \w{1} + \eta \w{2} - \e{1} \w{3} \right), \\
\dot{\e{}}_3 &=& \phantom{-} \frac{1}{2} \left(- \e{2} \w{1} + \e{1} \w{2} + \eta \w{3} \right) & , & 
\dot{\eta}    &=& -\frac{1}{2} \left( \e{1}\w{1} + \e{2}\w{2} + \e{3}\w{3} \right),
\end{array}
\end{equation}
where $\br{\C{B}}{\om{A}{B}} = [\w{1}~\w{2}~\w{3}]\tr$.  Equivalently, the angular velocity components are expressed in terms of the Euler parameters and their rates of change as
\begin{equation}\label{eq:QuatRate-2-Om}
\begin{array}{lcl}
\w{1} &=& 2 \left( \eta \dot{\e{}}_1 - \dot{\eta}\e{1} + \e{3}\dot{\e{}}_2 - \dot{\e{}}_3\e{2} \right), \\
\w{2} &=& 2 \left( \eta \dot{\e{}}_2 - \dot{\eta}\e{2} - \e{3}\dot{\e{}}_1 + \dot{\e{}}_3\e{1} \right), \\
\w{3} &=& 2 \left( \eta \dot{\e{}}_3 - \dot{\eta}\e{3} + \e{2}\dot{\e{}}_1 - \dot{\e{}}_2\e{1} \right).
\end{array}
\end{equation}
The relationships in Eqs.~\eqref{eq:C_BA-Quat}--\eqref{eq:QuatRate-2-Om} will be employed frequently in the derivation of Section~\ref{sect:Derivation}.

% -----------------------------------------------------------------
\section{Derivation of the Equations of Motion \label{sect:Derivation}}
Consider a particle of mass $m$ located at a point $P$ which moves along with the particle.  Let $\m{r}$ denote the position of $P$ measured relative to an inertially fixed point $O$.  Suppose $O$ is also fixed in frames $\C{E}$ and $\C{A}$ and point $P$ is fixed in frame $\C{B}$.  Here, frame $\C{E}$ is the observation frame (central, rotating body frame) in which the relative velocity, $^\C{E}{\m{v}}$, is sought.  In addition, the position and velocity frames, $\C{A}$ and $\C{B}$, are chosen such that the basis vectors $\h{a}_1$ and $\h{b}_1$ are aligned with $\m{r}$ and $^\C{E}{\m{v}}$ respectively.  The remaining degree of freedom for frame $\C{A}$ (rotation of $\h{a}_2$ and $\h{a}_3$ about $\m{r}$) and the remaining degree of freedom for frame $\C{B}$ (rotation of $\h{b}_2$ and $\h{b}_3$ about $^{\C{E}}\m{v}$) are removed by choosing any admissible initial orientation for the remaining basis vectors and applying constraints on $\w{A1}$ and $\w{B1}$, where $\br{A}{\om{E}{A}} = [\w{A1}~\w{A2}~\w{A3}]\tr$ and $\br{B}{\om{A}{B}} = [\w{B1}~\w{B2}~\w{B3}]\tr$.  Three sets of constraints on $\w{A1}$ and $\w{B1}$ are considered in Sections~\ref{sect:rv}--\ref{sect:rvh}, each producing distinct equations of motion with distinct properties.  The general case where $\w{A1}$ and $\w{B1}$ are arbitrary is considered in Section~\ref{sect:general} and forms the basis for the development of the parameterizations in Section~\ref{sect:guage}.

Throughout Sections~\ref{sect:general}--\ref{sect:guage} the following definitions are assumed.  Let $\{ \e{A1},\e{A2},\e{A3},\eta_A \}$ and $\{ \e{B1},\e{B2},\e{B3},\eta_B \}$ denote two sets of Euler parameters, termed here as $\CAE$ Euler parameters and $\CBA$ Euler parameters respectively.  The Euler parameters define the direction cosine matrices
\begin{equation}\label{eq:CAE-final}
\CAE = 
\begin{bmatrix} 
1 - 2(\e{A2}^2 + \e{A3}^2) & 
2(\e{A1} \e{A2} + \e{A3} \eta_A) &
2(\e{A1} \e{A3} - \e{A2} \eta_A) \\
2(\e{A2} \e{A1} - \e{A3} \eta_A) &
1 - 2(\e{A3}^2 + \e{A1}^2) &
2(\e{A2} \e{A3} + \e{A1} \eta_A) \\
2(\e{A3} \e{A1} + \e{A2} \eta_A) &
2(\e{A3} \e{A2} - \e{A1} \eta_A) &
1 - 2(\e{A1}^2 + \e{A2}^2)
\end{bmatrix}, \vspace{2mm}\\
\end{equation}
and
\begin{equation}\label{eq:CBA-final}
\CBA = 
\begin{bmatrix} 
1 - 2(\e{B2}^2 + \e{B3}^2) & 
2(\e{B1} \e{B2} + \e{B3} \eta_B) &
2(\e{B1} \e{B3} - \e{B2} \eta_B) \\
2(\e{B2} \e{B1} - \e{B3} \eta_B) &
1 - 2(\e{B3}^2 + \e{B1}^2) &
2(\e{B2} \e{B3} + \e{B1} \eta_B) \\
2(\e{B3} \e{B1} + \e{B2} \eta_B) &
2(\e{B3} \e{B2} - \e{B1} \eta_B) &
1 - 2(\e{B1}^2 + \e{B2}^2)
\end{bmatrix},
\end{equation}
where frames $\C{E}$, $\C{A}$, and $\C{B}$ are the aforementioned observation, position, and velocity frames.  The $\CAE$ Euler parameters, $\CBA$ Euler parameters, and the magnitudes $r$ and $v$ of vectors $\m{r}$ and $^\C{E}{\m{v}}$ comprise a ten parameter set that define the position and relative velocity of point $P$ relative to point $O$.  The differential equations of motion describing the evolution of the parameters are derived next.

\subsection{General Form\label{sect:general}}

The kinematic and kinetic equations are derived next by assuming $\w{A1}$ and $\w{B1}$ are defined by two arbitrary constraints.  Specific choices for these constraints are discussed in Sections~\ref{sect:rv}--\ref{sect:rvh}.  Throughout the derivation, it is assumed that the forces acting on the particle are thrust, lift, drag, and gravity.

\subsubsection{Kinematic Equations\label{sect:general-kinematics}}
The position and velocity of point $P$ relative to point $O$ as viewed by an observer in frame $\C{E}$ is described by
\begin{equation}\label{eq:rv}
\begin{array}{rcl}
\m{r} &=& r \h{a}_1, \\
^\C{E}\m{v} &=& v \h{b}_1,
\end{array}
\end{equation}
where $r$ is the magnitude of the position and $v$ is the speed of the particle as viewed by an observer in the observation frame $\C{E}$.  The relative velocity, $^\C{E}\m{v}$, is also defined by
\begin{equation}\label{eq:kinematics-1}
^\C{E}\m{v} = {\ddtm{\C{E}}{\m{r}}} = {\ddtm{\C{A}}{\m{r}}} + {\om{E}{A}} \times \m{r},
\end{equation}
where $\br{\C{A}}{\om{E}{A}} = [\w{A1}~\w{A2}~\w{A3}]\tr$.  Expressing Eq.~\eqref{eq:kinematics-1} in the basis fixed in frame $\C{A}$ yields
\begin{equation}\label{eq:kinematics-2}
\CAB \br{\C{B}}{^\C{E}\m{v}} = \br{\C{A}}{\ddtm{\C{A}}{\m{r}}}  + \brsk{\C{A}}{\om{E}{A}}  \{ \m{r} \}_{\C{A}},
\end{equation}
which is rewritten in matrix form and simplified to produce
\begin{equation}\label{eq:kinematics-final-0}
\begin{bmatrix}
\dot{r} \\ \phantom{-} r \w{A3} \\ -r \w{A2}
\end{bmatrix}
=
\CAB
\begin{bmatrix}
v \\ 0 \\ 0
\end{bmatrix}.
\end{equation}
Applying the definition of $\CAB$ and solving for $\dot{r}$, $\w{A2}$, and $\w{A3}$ yields the system of equations
\begin{equation}\label{eq:kinematics-final-1}
\begin{array}{lclclclclcl}
\dot{r}                   &=& v \left( 1 - 2(\e{B2}^2 + \e{B3}^2) \right) & , & 
\w{A2} &=& \frac{2 v}{r} \left( \eta_B \e{B2} - \e{B1} \e{B3} \right) & , & 
\w{A3} &=& \frac{2 v}{r} \left( \eta_B \e{B3} + \e{B1} \e{B2} \right), 
\end{array}
\end{equation}
where $\w{A1}$ is unconstrained in Eqs.~\eqref{eq:kinematics-final-0} and \eqref{eq:kinematics-final-1}.  Moreover, the ability to choose $\w{A1}$ is a direct consequence of the remaining degree of freedom in defining frame $\C{A}$ (rotations of $\{ \h{a}_2, \h{a}_3 \}$ about $\m{r}$) as discussed earlier.  Thus, a family of parameterizations exists, each with distinct properties, simply constraining $\w{A1}$ in different ways.  For example, suppose $\w{A1}$ is arbitrary.  The Euler parameter rates are then determined by Eq.~\eqref{eq:QuatRate} as
\begin{equation}\label{eq:kinematics-final-2}
\begin{array}{lclclcl}
\dot{\e{}}_{A1} &=& \phantom{-} \frac{1}{2} \left( \phantom{-} \eta_A \w{A1} - \e{A3} \w{A2} + \e{A2} \w{A3} \right) & , & 
\dot{\e{}}_{A2} &=& \phantom{-} \frac{1}{2} \left( \e{A3} \w{A1} + \eta_A \w{A2} - \e{A1} \w{A3} \right), \\
\dot{\e{}}_{A3} &=& \phantom{-} \frac{1}{2} \left(- \e{A2} \w{A1} + \e{A1} \w{A2} + \eta_A \w{A3} \right) & , & 
\dot{\eta}_A    &=& -\frac{1}{2} \left( \e{A1}\w{A1} + \e{A2}\w{A2} + \e{A3}\w{A3} \right),
\end{array}
\end{equation}
Equations~\eqref{eq:kinematics-final-1} and \eqref{eq:kinematics-final-2} then define the kinematics.  Although five differential equations define the motion of the parameters $\{ r, \e{A1}, \e{A2}, \e{A3}, \eta_{A} \}$, one degree of freedom is removed via the constraint $\e{A1}^2 + \e{A2}^2 + \e{A3}^2 + \eta_A^2 = 1$ (which is implicit in Eq.~\eqref{eq:kinematics-final-2}), while another degree of freedom is removed with an appropriate constraint on $\w{A1}$.  

\subsubsection{Kinetic Equations\label{sect:general-kinetics}}

Consider now the kinetics.  Newton's second law is then expanded in terms of the relative velocity as
\begin{equation}\label{eq:dynamics-1}
\frac{1}{m}\m{F} = {\ddtm{\C{B}}{^\C{E}\m{v}}} 
+ \left( \om{A}{B} + \om{E}{A} + 2 {\om{N}{E}} \right) \times {^\C{E}\m{v}}
+ \al{N}{E} \times \m{r}
+ \om{N}{E} \times \om{N}{E} \times \m{r},
\end{equation}
where $\br{\C{B}}{\om{A}{B}} = [\w{B1}~\w{B2}~\w{B3}]\tr$, $\br{\C{A}}{\om{E}{A}} = [\w{A1}~\w{A2}~\w{A3}]\tr$, $\br{\C{E}}{\om{N}{E}} = [0~0~\w{e}]\tr$, and $\al{N}{E} = \m{0}$.  Note that the previous statement assumes the central body rotates with constant angular velocity $\w{e}$ about $\h{e}_3$.  

\begin{figure}[hbt!]
  \centering
  \begin{tabular}{lr}
  \subfloat[Rear view ($\alpha = \delta = 0$).]{\includegraphics[width=.475\textwidth]{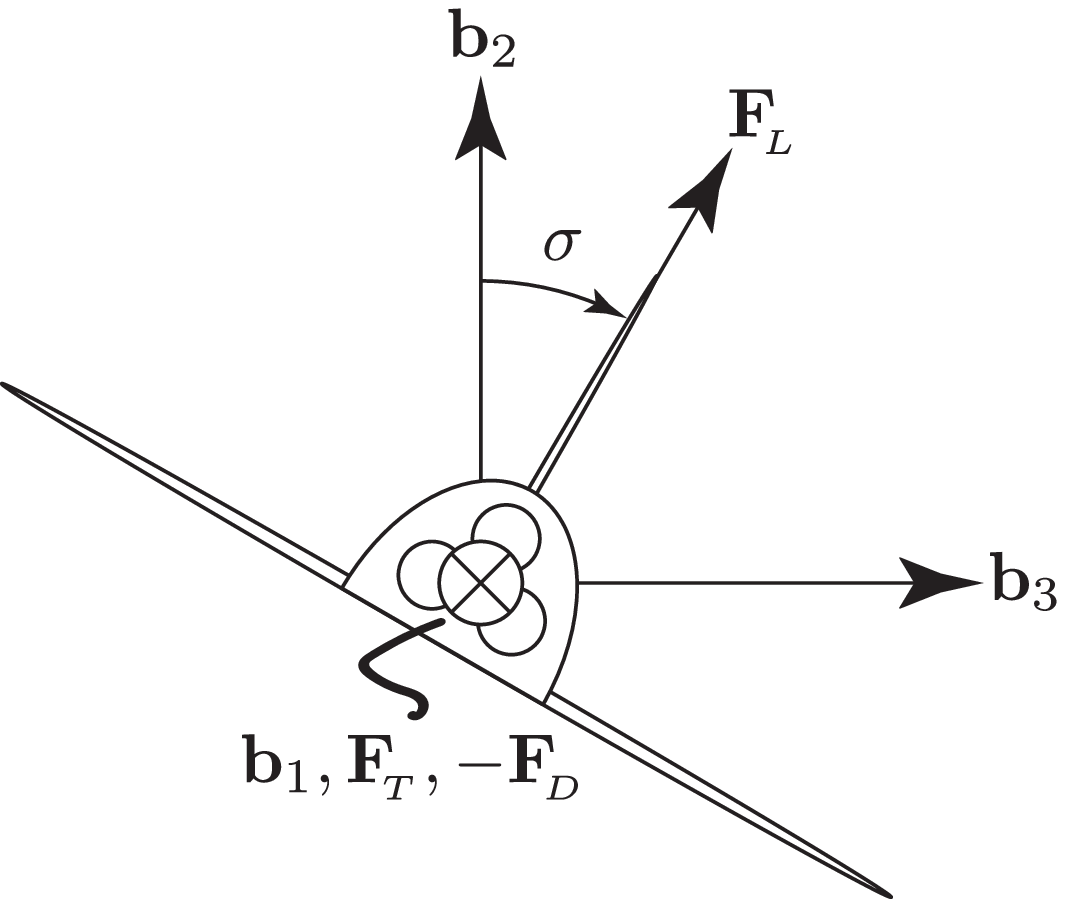}}
  &
  \subfloat[Side view ($\sigma = 0$).]{\includegraphics[width=.475\textwidth]{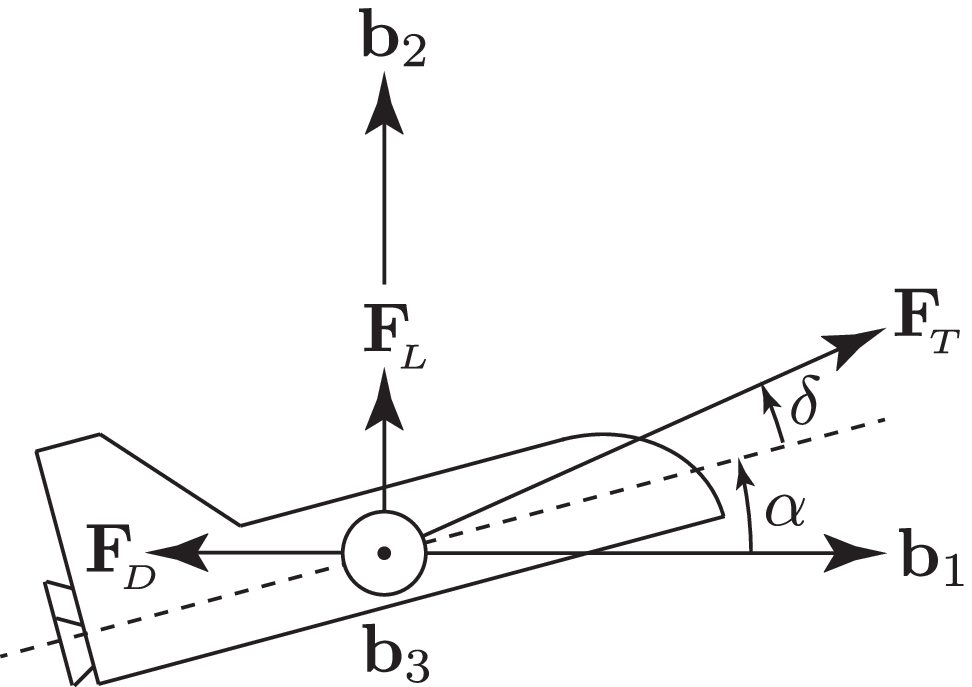}}
  \end{tabular}
  \caption{Free body diagram (excluding gravity) and assuming $L$ is positive.}\label{fig:FBD}
\end{figure}

Next, suppose that thrust, lift, drag, and gravity are the forces acting on the particle, denoted $\m{F}_T$, $\m{F}_L$, $\m{F}_D$, and $\m{F}_g$ respectively.  Suppose further that the thrust vector lies in the lift-drag plane and that the atmosphere is fixed in frame $\C{E}$ (fixed to the central, rotating body).  The net force acting on the particle is then expressed as
\begin{equation}\label{eq:netForce}
\m{F} = \left( T \cos (\alpha + \delta) - D \right) \h{b}_1
          + \left( T \sin (\alpha + \delta) +L \right) \cos (\sigma) \h{b}_2
          + \left( T \sin (\alpha + \delta) +L \right) \sin (\sigma) \h{b}_3
          - \frac{\mu_e}{r^2} \h{a}_1,
\end{equation}
where $T$, $L$, and $D$ are the thrust, lift, and drag force magnitudes (noting that $L$ may be positive or negative), $\alpha$ is the angle of attack, $\sigma$ is the bank angle (measured as the angle of rotation about $\h{b}_1$ from $\h{b}_2$ to the positive lift direction), $\delta$ is the offset of the thrust vector from the body x-axis, and $\mu_e$ is the gravitational parameter of the central, rotating body.  Figure~\ref{fig:FBD} illustrates the free body diagram (excluding gravity).

Next, the net force, denoted $\tilde{\m{F}}$, is related to the resultant force acting on the particle via the relationship 
\begin{equation}\label{eq:F-apparent-1}
\tilde{\m{F}} = \m{F} 
- m \left(
2 {\om{N}{E}} \times {^\C{E}\m{v}}
+ \om{N}{E} \times \om{N}{E} \times \m{r}
\right),
\end{equation}
such that $\tilde{\m{F}}/m = {^\C{E}{\m{a}}}$ and it has been assumed that $\al{N}{E} = \m{0}$.  Expressing the result of Eq.~\eqref{eq:F-apparent-1} in the basis fixed in frame $\C{B}$ and dividing by $m$ yields
\begin{equation}\label{eq:F-apparent-2}
\begin{array}{rcl}
\frac{1}{m}\br{\C{B}}{\tilde{\m{F}}} &\hspace{-2mm} =& \hspace{-2mm} \frac{1}{m} \br{\C{B}}{\m{F}}
- 2 \CBE \brsk{\C{E}}{\om{N}{E}} \CEB \br{\C{B}}{^\C{E}\m{v}} 
- \CBE \brsk{\C{E}}{\om{N}{E}} \brsk{\C{E}}{\om{N}{E}} \CEA \br{\C{A}}{\m{r}}. 
\end{array}
\end{equation}
Expressing this last in matrix form and simplifying gives
\begin{equation}\label{eq:F-apparent-3}
\begin{array}{rcl}
\frac{1}{m} \begin{bmatrix}
\tilde{f}_1 \\ \tilde{f}_2 \\ \tilde{f}_3
\end{bmatrix}
&=&
\frac{1}{m}\begin{bmatrix}
f_1 \\ f_2 \\ f_3
\end{bmatrix}

- 2 \w{e} v
\begin{bmatrix}
0 \\
\phantom{-}\CBE(3,3) \\
-\CBE(2,3)
\end{bmatrix}

- r \w{e}^2 \CBA
\begin{bmatrix}
\left( \CAE(1,3) \right)^2  - 1 \\
\CAE(1,3) \CAE(2,3) \\
\CAE(1,3) \CAE(3,3)
\end{bmatrix},
\end{array}
\end{equation}
where $\br{\C{B}}{\tilde{\m{F}}} = [\tilde{f}_1~\tilde{f}_2~\tilde{f}_3]\tr$, $\br{\C{B}}{{\m{F}}} = [{f}_1~{f}_2~{f}_3]\tr$, and $\CBE = \CBA \CAE$.  Notice that the force components $f_1$, $f_2$, and $f_3$ are given by
\begin{equation}\label{eq:netForce-B}
\begin{array}{rcl}
 \begin{bmatrix}
f_1 \\ f_2 \\ f_3
\end{bmatrix}

& = &
 
\begin{bmatrix}
T \cos (\alpha + \delta) - D \\ 
\left( T \sin (\alpha + \delta) +L \right) \cos (\sigma) \\
\left( T \sin (\alpha + \delta) +L \right) \sin (\sigma)
\end{bmatrix}
\displaystyle -\frac{\mu_e}{r^2} \CBA(:,1).
\end{array}
\end{equation}

Now return to Eq.~\eqref{eq:dynamics-1}.  Substituting the definition for $\tilde{\m{F}}$, expressing all vectors in the basis fixed in frame $\C{B}$, and reordering terms produces
\begin{equation}\label{eq:dynamics-2}
\br{\C{B}}{\ddtm{\C{B}}{^\C{E}\m{v}}} 
+ \brsk{\C{B}}{\om{A}{B}} \br{\C{B}}{^\C{E}\m{v}}
=
\frac{1}{m} \br{\C{B}}{\tilde{\m{F}}}
- \CBA \brsk{\C{A}}{\om{E}{A}} \CAB \br{\C{B}}{^\C{E}\m{v}},
\end{equation}
which is equivalent to
\begin{equation}\label{eq:dynamics-final-0}
\begin{array}{rcl}
\begin{bmatrix}
\dot{v} \\ \phantom{-} v \w{B3} \\ -v \w{B2}
\end{bmatrix}
&=&
\frac{1}{m} \begin{bmatrix}
\tilde{f}_1 \\ \tilde{f}_2 \\ \tilde{f}_3
\end{bmatrix}

- v
\begin{bmatrix}
\m{0}\tr \\
\phantom{-}\CBA(3,:) \\
-\CBA(2,:)
\end{bmatrix}
\begin{bmatrix}
\w{A1} \\ \w{A2} \\ \w{A3}
\end{bmatrix}.
\end{array}
\end{equation}
Equation~\eqref{eq:dynamics-final-0} is solved for $\dot{v}$, $\w{B2}$, and $\w{B3}$ to produce
\begin{equation}\label{eq:dynamics-final-1}
\begin{array}{rcl}
\dot{v} &=& \displaystyle \frac{1}{m} \tilde{f}_1, \\
\w{B2} &=& \displaystyle - \frac{1}{m v} \tilde{f}_3 
- \w{A1} \CBA(2,1)
- \frac{v}{r} \CBA(3,1), \\

\w{B3} &=& \displaystyle \frac{1}{m v} \tilde{f}_2
- \w{A1} \CBA(3,1)
+ \frac{v}{r} \CBA(2,1),
\end{array}
\end{equation}
%- \frac{v}{r} \CBA(3,1) = - \w{A2} \CBA(2,2) - \w{A3} \CBA(2,3),
%   \frac{v}{r} \CBA(2,1) = - \w{A2} \CBA(3,2) - \w{A3} \CBA(3,3),
where the definitions of $\w{A2}$ and $\w{A3}$ in Eq.~\eqref{eq:kinematics-final-1} have been used.
An equivalent expression of Eq.~\eqref{eq:dynamics-final-1} could be obtained solely in terms of the Euler parameters, but is omitted here for clarity and brevity.  Observing Eqs.~\eqref{eq:dynamics-final-0} and \eqref{eq:dynamics-final-1} it is noticed that $\w{B1}$ is unconstrained.  The freedom to choose $\w{B1}$ stems from the remaining degree of freedom in the definition of frame $\C{B}$ (That is, rotations of $\{ \h{b}_2, \h{b}_3 \}$ about $^{\C{E}}\m{v}$ are arbitrary).  Following the same reasoning as in Section~\ref{sect:general-kinematics}, suppose $\w{B1}$ is arbitrary for now.  The Euler parameter rates are then determined by Eq.~\eqref{eq:QuatRate} as
\begin{equation}\label{eq:dynamics-final-2}
\begin{array}{lclclcl}
\dot{\e{}}_{B1} &=& \phantom{-} \frac{1}{2} \left( \phantom{-} \eta_B \w{B1} - \e{B3} \w{B2} + \e{B2} \w{B3} \right) & , & 
\dot{\e{}}_{B2} &=& \phantom{-} \frac{1}{2} \left( \e{B3} \w{B1} + \eta_B \w{B2} - \e{B1} \w{B3} \right), \\
\dot{\e{}}_{B3} &=& \phantom{-} \frac{1}{2} \left(- \e{B2} \w{B1} + \e{B1} \w{B2} + \eta_B \w{B3} \right) & , & 
\dot{\eta}_B    &=& -\frac{1}{2} \left( \e{B1}\w{B1} + \e{B2}\w{B2} + \e{B3}\w{B3} \right).
\end{array}
\end{equation}
Together, Eqs.~\eqref{eq:dynamics-final-1}, and \eqref{eq:dynamics-final-2} define the kinetic equations.  The kinetic equations define the evolution of five parameters, $\{ v, \e{B1}, \e{B2}, \e{B3}, \eta_{B} \}$, in addition to the five introduced in the kinematics derivation, $\{ r, \e{A1}, \e{A2}, \e{A3}, \eta_{A} \}$.  Once again, it is noted that there are only three degrees of freedom, because one degree of freedom is removed by the constraint $\e{B1}^2 + \e{B2}^2 + \e{B3}^2 + \eta_B^2 = 1$ (implicit in Eq.~\eqref{eq:dynamics-2})and one more degree of freedom is removed once a constraint on $\w{B1}$ is applied.

\subsubsection{Bank Angle Conversion\label{sect:general-bank}}

\begin{figure}[hbt!]
  \centering
  \includegraphics[width=.3\textwidth]{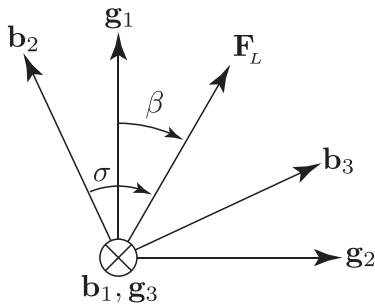}
  \caption{Relationship between $\beta$ and $\sigma$ (positive $L$ assumed).}\label{fig:bankBS}
\end{figure}

Section~\ref{sect:general-kinetics} defines the bank angle as the angle of rotation about $\h{b}_1$ from $\h{b}_2$ to the positive lift direction.  The previous definition remains well defined even in vertical flight ($\h{a}_1$ and $\h{b}_1$ co-linear).  In contrast, traditional definitions of the bank angle rely upon the $\{ \m{r}, ^{\C{E}}\m{v} \}$ plane for reference when measuring the bank angle.  For instance, the following definition is commonly used.  Let $\beta$ denote the bank angle measured as the angle of rotation about $\h{g}_3$ from $\h{g}_1$ to the positive lift direction, where the orthonormal basis $\{ \h{g}_1, \h{g}_2, \h{g}_3 \}$ is defined as
\begin{equation}\label{eq:gBasis}
\begin{array}{lclclclclcl}
\h{g}_1   &=& \h{g}_2  \times \h{g}_3  & , & 
\h{g}_2  &=& \displaystyle - \frac{\m{r} \times ^{\C{E}}\m{v}}{|| \m{r} \times ^{\C{E}}\m{v} ||} & , & 
\h{g}_3  &=& \displaystyle \frac{^{\C{E}}\m{v}}{|| ^{\C{E}}\m{v} ||}.
\end{array}
\end{equation}
Figure~\ref{fig:bankBS} provides a visual of the relationship between $\sigma$ and $\beta$.

While the definition for $\beta$ is not well defined in vertical flight, it is more physically meaningful than the definition for $\sigma$.  Thus, it is useful (and insightful) to derive a mapping from $\sigma$ to $\beta$.  The map is derived as follows.  Begin by noticing that $\h{g}_3 = \h{b}_1$ and that $\h{g}_2$ may be expressed as
\begin{equation}\label{eq:g2-1}
\h{g}_2 = -\frac{\h{a}_1 \times \h{b}_1}{|| \h{a}_1 \times \h{b}_1 ||}.
\end{equation}
Expressing all vectors in Eq.~\eqref{eq:g2-1} in the basis fixed in frame $\C{B}$ and simplifying yields
\begin{equation}\label{eq:g2-2}
\br{\C{B}}{\h{g}_2} = \frac{1}{\sqrt{1 - \left( \CBA(1,1) \right)^2}} \begin{bmatrix}
0 \\ -\CBA(3,1) \\ \phantom{-}\CBA(2,1)
\end{bmatrix}.
\end{equation}
Likewise, $\h{g}_1$ is expressed in the basis fixed in frame $\C{B}$ as
\begin{equation}\label{eq:g1}
\br{\C{B}}{\h{g}_1} = \frac{1}{\sqrt{1 - \left( \CBA(1,1) \right)^2}} \begin{bmatrix}
0 \\ \CBA(2,1) \\ \CBA(3,1)
\end{bmatrix}.
\end{equation}
Next, observe Fig.~\ref{fig:bankBS} and notice that the positive lift direction may be equivalently expressed as $\m{F}_L / L = \cos (\beta) \h{g}_1 + \sin (\beta) \h{g}_2$ or as $\m{F}_L / L = \cos (\sigma) \h{b}_2 + \sin (\sigma) \h{b}_3$.  It follows that
\begin{equation}\label{eq:BS-1}
\beta = \textrm{atan2}\left( \frac{\m{F}_L}{L} \cdot \h{g}_2~~,~~\frac{\m{F}_L}{L} \cdot \h{g}_1 \right),
\end{equation}
where atan2($y$,$x$) is the four-quadrant inverse tangent operator.  Expressing the vectors in Eq.~\eqref{eq:BS-1} in the basis fixed in frame $\C{B}$ and simplifying yields the relationship
\begin{equation}\label{eq:BS-2}
\beta = \textrm{atan2}\left( \sin(\sigma) \CBA(2,1) - \cos(\sigma) \CBA(3,1)~~,~~\cos(\sigma) \CBA(2,1) + \sin(\sigma) \CBA(3,1) \right).
\end{equation}

The map from $\sigma$ to $\beta$ in Eq.~\eqref{eq:BS-2} may be differentiated to provide the relationship between the bank angle rates $\dot{\sigma}$ and $\dot{\beta}$ as well.  Differentiating Eq.~\eqref{eq:BS-2} and simplifying extensively produces the relationship
\begin{equation}\label{eq:BSdot-1}
\dot{\beta} = \left( \dot{\sigma} + \w{B1} \right) - \frac{\CBA(1,1)}{1 - \left( \CBA(1,1) \right)^2} \left( \w{B2} \CBA(2,1) + \w{B3} \CBA(3,1) \right),
\end{equation}
where it is noted that the relationship $\dot{\mathbf{C}}_{\C{B} \C{A}} = -\brsk{\C{B}}{\om{\C{A}}{\C{B}}} \CBA$ has been used.  
%Finally, it is noted that Eq.~\eqref{eq:BSdot-1} may also be written in a more physically meaningful form as
%\begin{equation}\label{eq:BSdot-1}
%\dot{\beta} = \left( \dot{\sigma} + \w{B1} \right) - \frac{\h{a}_1 \cdot \h{b}_1}{1 - \left( \h{a}_1 \cdot \h{b}_1 \right)^2} \left( \om{\C{A}}{\C{B}} \cdot \left( \h{a}_1 - (\h{a}_1 \cdot \h{b}_1)\h{b}_1 \right) \right).
%\end{equation}

\subsection{Angular Velocity Constraints\label{sect:guage}}

In Section~\ref{sect:general} it was found that $\w{A1}$ and $\w{B1}$ could be arbitrarily chosen.  While it may seem a trivial task to choose two constraints that define $\w{A1}$ and $\w{B1}$, the nuances of such choices have major implications for the equations of motion (as will become clear shortly).  Three potentially useful combinations of the constraints on $\w{A1}$ and $\w{B1}$ are now explored.  Each pair of constraints gives rise to a distinct parameterization with distinct properties.  For ease of communication, the three parameterizations are termed here as the $rv$-, $rvL$-, and $rvh$-Euler parameterizations.

\subsubsection{The $rv$-Euler Parameterization\label{sect:rv}}
The $rv$-Euler parameterization applies perhaps the simplest and most obvious constraints on $\w{A1}$ and $\w{B1}$, setting them both equal to zero.  The constraints are equivalent to
\begin{equation}\label{eq:wcon-rv}
\begin{array}{rcl}
0 &=& \eta_A \edot{A1} - \dot{\eta}_A \e{A1} + \e{A3} \edot{A2} - \edot{A3} \e{A2}, \\
0 &=& \eta_B \edot{B1} - \dot{\eta}_B \e{B1} + \e{B3} \edot{B2} - \edot{B3} \e{B2},
\end{array}
\end{equation}
by Eq.~\eqref{eq:QuatRate-2-Om} and are implicitly applied in the equations of motion by removing all $\w{A1}$ and $\w{B1}$ terms from Eqs.~\eqref{eq:kinematics-final-2}, \eqref{eq:dynamics-final-1}, and \eqref{eq:dynamics-final-2} in Section~\ref{sect:general}.  The kinematic equations are then summarized by
\begin{equation}\label{eq:rv-kinematics1}
\begin{array}{lclclclclcl}
\dot{r}                   &=& v \left( 1 - 2(\e{B2}^2 + \e{B3}^2) \right) & , & 
\w{A2} &=& \frac{2 v}{r} \left( \eta_B \e{B2} - \e{B1} \e{B3} \right) & , & 
\w{A3} &=& \frac{2 v}{r} \left( \eta_B \e{B3} + \e{B1} \e{B2} \right),
\end{array}
\end{equation}
and
\begin{equation}\label{eq:rv-kinematics2}
\begin{array}{lclclcl}
\dot{\e{}}_{A1} &=& \phantom{-} \frac{1}{2} \left(-\e{A3} \w{A2} + \e{A2} \w{A3} \right) & , & 
\dot{\e{}}_{A2} &=& \phantom{-} \frac{1}{2} \left( \eta_A \w{A2} - \e{A1} \w{A3} \right), \\
\dot{\e{}}_{A3} &=& \phantom{-} \frac{1}{2} \left( \phantom{-} \e{A1} \w{A2} + \eta_A \w{A3} \right) & , & 
\dot{\eta}_A    &=& -\frac{1}{2} \left( \e{A2}\w{A2} + \e{A3}\w{A3} \right).
\end{array}
\end{equation}
Likewise, the kinetic equations are summarized as
\begin{equation}\label{eq:rv-dynamics1}
\begin{array}{rcl}
\dot{v} &=& \frac{1}{m} \tilde{f}_1, \\

\w{B2} &=& - \frac{1}{m v} \tilde{f}_3 
-\frac{2v}{r} (\e{B1}\e{B3} + \e{B2}\eta_B), \\

\w{B3} &=& \frac{1}{m v} \tilde{f}_2
+ \frac{2v}{r} (\e{B1}\e{B2} - \e{B3}\eta_B),
\end{array}
\end{equation}
and
\begin{equation}\label{eq:rv-dynamics2}
\begin{array}{lclclcl}
\dot{\e{}}_{B1} &=& \phantom{-} \frac{1}{2} \left(-\e{B3} \w{B2} + \e{B2} \w{B3} \right) & , & 
\dot{\e{}}_{B2} &=& \phantom{-} \frac{1}{2} \left( \eta_B \w{B2} - \e{B1} \w{B3} \right), \\
\dot{\e{}}_{B3} &=& \phantom{-} \frac{1}{2} \left(\phantom{-} \e{B1} \w{B2} + \eta_B \w{B3} \right) & , & 
\dot{\eta}_B    &=& -\frac{1}{2} \left( \e{B2}\w{B2} + \e{B3}\w{B3} \right),
\end{array}
\end{equation}
where $\tilde{f}_1$, $\tilde{f}_2$, and $\tilde{f}_3$ remain defined by Eq.~\eqref{eq:F-apparent-3}.  Finally, it is noted that the bank angle map of Eq.~\eqref{eq:BS-2} remains unchanged and that the bank angle rate relationship in Eq.~\eqref{eq:BSdot-1} simplifies to
\begin{equation}\label{eq:rv-BSdot}
\dot{\beta} = \dot{\sigma} - \frac{\CBA(1,1)}{1 - \left( \CBA(1,1) \right)^2} \left( \w{B2} \CBA(2,1) + \w{B3} \CBA(3,1) \right),
\end{equation}

\subsubsection{The $rvL$-Euler Parameterization\label{sect:rvL}}
The $rvL$-Euler parameterization takes a more nuanced approach.  Suppose $\w{A1} = 0$ once again.  Thus, the kinematic relations of Eqs.~\eqref{eq:rv-kinematics1} and \eqref{eq:rv-kinematics2} still hold.  Next, suppose $\w{B1}$ is a control variable and suppose the basis vector $\h{b}_2$ defines the positive lift direction.  The net force components in Eq.~\eqref{eq:netForce-B} are then expressed more compactly as
\begin{equation}\label{eq:rvL-netForce}
\begin{array}{rcl}
 \begin{bmatrix}
f_1 \\ f_2 \\ f_3
\end{bmatrix}

& = &
 
\begin{bmatrix}
T \cos (\alpha + \delta) - D \\ 
\left( T \sin (\alpha + \delta) +L \right) \\
0
\end{bmatrix}

-\frac{\mu_e}{r^2} \CBA(:,1),
\end{array}
\end{equation}
where it is noted that the bank angle $\sigma$ no longer appears.  Given $\w{A1} = 0$ and $\w{B1}$ is a control variable, the kinetic equations are summarized as
\begin{equation}\label{eq:rvL-dynamics1}
\begin{array}{rcl}
\dot{v} &=& \frac{1}{m} \tilde{f}_1, \\

\w{B2} &=& - \frac{1}{m v} \tilde{f}_3 
- \frac{2v}{r} (\e{B1}\e{B3} + \e{B2}\eta_B), \\

\w{B3} &=& \frac{1}{m v} \tilde{f}_2
+ \frac{2v}{r} (\e{B1}\e{B2} - \e{B3}\eta_B),
\end{array}
\end{equation}
and
\begin{equation}\label{eq:rvL-dynamics2}
\begin{array}{lclclcl}
\dot{\e{}}_{B1} &=& \phantom{-} \frac{1}{2} \left( \phantom{-} \eta_B \w{B1} - \e{B3} \w{B2} + \e{B2} \w{B3} \right) & , & 
\dot{\e{}}_{B2} &=& \phantom{-} \frac{1}{2} \left( \e{B3} \w{B1} + \eta_B \w{B2} - \e{B1} \w{B3} \right), \\
\dot{\e{}}_{B3} &=& \phantom{-} \frac{1}{2} \left(- \e{B2} \w{B1} + \e{B1} \w{B2} + \eta_B \w{B3} \right) & , & 
\dot{\eta}_B    &=& -\frac{1}{2} \left( \e{B1}\w{B1} + \e{B2}\w{B2} + \e{B3}\w{B3} \right),
\end{array}
\end{equation}
where $\tilde{f}_1$, $\tilde{f}_2$, and $\tilde{f}_3$ remain defined by Eq.~\eqref{eq:F-apparent-3}.  Notice that the bank angle $\sigma$ is entirely removed from the equations of motion.  Finally, it is noted that the bank angle and bank angle rate relationships of Eqs.~\eqref{eq:BS-2} and \eqref{eq:BSdot-1} simplify to
\begin{equation}\label{eq:rvL-BS}
\beta = \textrm{atan2}\left( -\CBA(3,1)~~,~~ \CBA(2,1) \right),
\end{equation}
and
\begin{equation}\label{eq:rvL-BSdot}
\dot{\beta} = \w{B1} - \frac{\CBA(1,1)}{1 - \left( \CBA(1,1) \right)^2} \left( \w{B2} \CBA(2,1) + \w{B3} \CBA(3,1) \right),
\end{equation}
respectively.  Notice that the bank angle $\sigma$ does not appear and that $\w{B1}$ is effectively a bank angle rate command.

\subsubsection{The $rvh$-Euler Parameterization\label{sect:rvh}}
The $rvh$-Euler parameterization offers a simpler set of equations, but at a cost.  The equations of motion contain a singularity in vertical flight.  While still useful in non-vertical flight applications, the $rvh$-Euler parameterization also serves as a reminder that the constraints placed on $\w{A1}$ and $\w{B1}$ should be chosen with care.  The constraints on $\w{A1}$ and $\w{B1}$ in the $rvh$-Euler parameterization are derived as follows.

Suppose both the $\h{a}_3$ and $\h{b}_3$ basis vectors are directed along the specific angular momentum $^\C{E}{\m{h}} = \m{r} \times ^\C{E}{\m{v}}$.  Notice that the previous statement implies that the bases $\{ \h{a}_1, \h{a}_2, \h{a}_3 \}$ and $\{ \h{b}_1, \h{b}_2, \h{b}_3 \}$ are offset from one another by a simple rotation about $\h{a}_3 = \h{b}_3$.  Thus, the definition of the direction cosine matrix $\CBA$ in Eq.~\eqref{eq:CBA-final} simplifies to
\begin{equation}\label{eq:rvh-CBA}
\CBA = 
\begin{bmatrix} 
1 - 2 \e{B3}^2 & 
2 \e{B3} \eta_B &
0 \\
-2 \e{B3} \eta_B &
1 - 2 \e{B3}^2 &
0 \\
0 &
0 &
1
\end{bmatrix},
\end{equation}
where it is noted that $\e{B1}$ and $\e{B2}$ are identically zero.  Given that $\h{a}_3$ and $^\C{E}{\m{h}}$ are aligned at all times, the condition
\begin{equation}\label{eq:rvh-hcon-A-0}
\ddtm{\C{E}}{\frac{^\C{E}{\m{h}}}{h}} = \ddtm{\C{E}}{\h{a}_3},
\end{equation}
must hold, noting that $h$ is the magnitude of $^\C{E}{\m{h}}$.  Expanding both sides of Eq.~\eqref{eq:rvh-hcon-A-0} produces
\begin{equation}\label{eq:rvh-hcon-A-1}
\left( \frac{1}{h} \right) \ddtm{\C{E}}{^\C{E}{\m{h}}}
- \left( \frac{^\C{E}{\m{h}}}{h^2} \right) \frac{d}{dt}(h)
=
\ddtm{\C{A}}{\h{a}_3} + \om{E}{A} \times \h{a}_3,
\end{equation}
where
\begin{equation}\label{eq:rvh-ddt(h)}
\begin{array}{rcl}
\ddt{\C{E}}{^\C{E}{\m{h}}}
&=&
\m{r} \times {^\C{E}{\m{a}}}, \\

\frac{d}{dt} \left( h \right)
&=&
\h{a}_3 \cdot {\ddt{\C{E}}{^\C{E}{\m{h}}}},
\end{array}
\end{equation}
and noting that $^\C{E}{\m{a}}$ is the acceleration of point $P$ relative to point $O$ as viewed by an observer in the observation frame $\C{E}$.  Substituting Eq.~\eqref{eq:rvh-ddt(h)} into Eq.~\eqref{eq:rvh-hcon-A-1} and further simplifying yields
\begin{equation}\label{eq:rvh-hcon-A-2}
\frac{1}{h} \left(
\m{r} \times {^\C{E}{\m{a}}}
- \h{a}_3 \left( \h{a}_3 \cdot \m{r} \times {^\C{E}{\m{a}}} \right) \right)
=
\w{A2} \h{a}_1 - \w{A1} \h{a}_2,
\end{equation}
where $\br{\C{A}}{\om{E}{A}} \equiv [\w{A1}~\w{A2}~\w{A3}]\tr$.  Next, recall that the apparent force, denoted $\m{\tilde{F}}$, satisfies $\m{\tilde{F}} = m~{^{\C{E}}{\m{a}}}$.  Substituting ${^{\C{E}}{\m{a}}} = \frac{1}{m} \m{\tilde{F}}$ and $\m{r} = r \h{a}_1$ into Eq.~\eqref{eq:rvh-hcon-A-2}, expressing all vectors in the basis fixed in frame $\C{A}$, and simplifying produces
\begin{equation}\label{eq:rvh-hcon-A-3}
\begin{bmatrix}
\w{A2} \\ \w{A1} \\ 0
\end{bmatrix}
=
\frac{r}{h m}
\begin{bmatrix}
0 \\ \tilde{f}_{3} \\ 0
\end{bmatrix},
\end{equation}
where $\br{\C{B}}{\m{\tilde{F}}} = [ \tilde{f}_1, \tilde{f}_2, \tilde{f}_3 ]\tr$ and noting that $\tilde{f}_3 \h{a}_3 = \tilde{f}_3 \h{b}_3$ because $\h{a}_3$ is aligned with $\h{b}_3$.  Equation~\eqref{eq:rvh-hcon-A-3} produces the conditions
\begin{equation}\label{eq:rvh-hcon-A-final-1}
\begin{array}{rcl}
\w{A1} &=& \displaystyle \frac{r}{h m} \tilde{f}_{3}, \\
\w{A2} &=& 0.
\end{array}
\end{equation}

Notice in Eq.~\eqref{eq:rvh-hcon-A-final-1} that the constraint on $\w{A1}$ includes the magnitude of the specific angular momentum.  An equivalent expression without $h$ is derived as follows.  Recall that $^\C{E}{\m{h}} = \m{r} \times ^\C{E}{\m{v}}$ and notice that $h$ can be expressed as
\begin{equation}
h = r v \sin \phi,
\end{equation}
where $\phi \in [0,\pi]$ is the angle between $\m{r}$ and $^{\C{E}}\m{v}$.  Notice also that a simple rotation about $\h{a}_3 = \h{b}_3$ by the same angle $\phi$ aligns the basis vectors fixed in frame $\C{A}$ with those fixed in frame $\C{B}$.  Thus, $\CBA$ is expressed in terms of the angle $\phi$ as
\begin{equation}\label{eq:rvh-CBA-phi}
\CBA = 
\begin{bmatrix} 
\phantom{-}\cos \phi & 
\sin \phi &
0 \\
-\sin \phi &
\cos \phi &
0 \\
\phantom{-}0 &
0 &
1
\end{bmatrix}.
\end{equation}
Comparing Eq.~\eqref{eq:rvh-CBA-phi} with Eq.~\eqref{eq:rvh-CBA} it is observed that $\sin \phi = 2 \e{B3} \eta_B$.  Therefore, an equivalent expression for $h$ is
\begin{equation}\label{eq:rvh-hmag-def}
h = 2 r v \e{B3} \eta_B,
\end{equation}
and the constraint on $\w{A1}$ in Eq.~\eqref{eq:rvh-hcon-A-final-1} is rewritten as
\begin{equation}\label{eq:rvh-hcon-wA1}
\w{A1} = \frac{\tilde{f}_{3}}{2 m v \e{B3} \eta_B}.
\end{equation}

Given the $\w{A1}$ constraint of Eq.~\eqref{eq:rvh-hcon-wA1} and the updated definition of $\CBA$ in Eq.~\eqref{eq:rvh-CBA}, the kinematic equations are summarized as
\begin{equation}\label{eq:rvh-kinematics1}
\begin{array}{rcl}
\dot{r}  &=& v \left( 1 - 2 \e{B3}^2 \right), \\
\w{A1} &=& \frac{\tilde{f}_{3}}{2 m v \e{B3} \eta_B}, \\
\w{A3} &=& \frac{2 v}{r} \eta_B \e{B3},
\end{array}
\end{equation}
and
\begin{equation}\label{eq:rvh-kinematics2}
\begin{array}{rcl}
\edot{A1} &=& \phantom{-} \frac{1}{2} \w{A1} \eta_A + \frac{1}{2} \w{A3} \e{A2}, \\
\edot{A2} &=& \phantom{-} \frac{1}{2} \w{A1} \e{A3} - \frac{1}{2} \w{A3} \e{A1}, \\
\edot{A3} &=&                  -  \frac{1}{2} \w{A1} \e{A2} + \frac{1}{2} \w{A3} \eta_A, \\
\dot{\eta}_A &=&              -  \frac{1}{2} \w{A1} \e{A1} - \frac{1}{2} \w{A3} \e{A3},
\end{array}
\end{equation}
noting that $\w{A2} = 0$ in Eqs.~\eqref{eq:rvh-kinematics1} and \eqref{eq:rvh-kinematics2}.  Next, the kinematic equations are summarized as
\begin{equation}\label{eq:rvh-dynamics1}
\begin{array}{rcl}
\dot{v} &=& \frac{1}{m} \tilde{f}_1, \\
\w{B3} &=& \frac{1}{m v} \tilde{f}_2 -\frac{2 v}{r} \eta_B \e{B3},
\end{array}
\end{equation}
and
\begin{equation}\label{eq:rvh-dynamics2}
\begin{array}{rcl}
\edot{B3} &=& \phantom{-} \frac{1}{2} \w{B3} \eta_B, \\
\dot{\eta}_B &=&              -   \frac{1}{2} \w{B3} \e{B3},
\end{array}
\end{equation}
where it is noted that $\w{B1} = 0$ and $\w{B2} = 0$ because the basis fixed in frame $\C{B}$ is offset from the basis fixed in frame $\C{A}$ by a simple rotation about $\h{a}_3 = \h{b}_3$ (that is, $\br{\C{B}}{\om{A}{B}} = [0~0~\w{B3}]\tr$).  It is also noted that the apparent force components $\tilde{f}_1$, $\tilde{f}_2$, and $\tilde{f}_3$ remain defined by Eq.~\eqref{eq:F-apparent-3}.  Finally, the bank angle and bank angle rate relations in Eqs.~\eqref{eq:BS-2} and \eqref{eq:BSdot-1} reduce to
\begin{equation}\label{eq:rvh-BS}
\begin{array}{rcl}
\beta &=& \sigma \pm \pi, \\
\dot{\beta} &=& \dot{\sigma}.
\end{array}
\end{equation}

% -----------------------------------------------------------------
\section{Example: Atmospheric Entry Trajectory Optimization \label{sect:Example}}
Consider the following variation of the atmospheric entry optimal control problem detailed in Refs.~\cite{Clarke2,Clarke1}.  The optimal control problems stated in Refs.~\cite{Clarke2,Clarke1} avoid enforcing an exact vertical impact condition (terminal flight path angle equal to $-90\deg$) due to singularity concerns.  Here, the $rv$-Euler parameterization is employed to allay singularity concerns and the vertical impact condition is enforced exactly.

\subsection{Problem Statement\label{sect:Ex-OCP}}
The atmospheric entry optimal control problem is stated in terms of the $rv$-Euler parameters as follows.  First, the objective functional to be minimized is given as
\begin{equation}\label{eq:Ex-Cost}
\C{J} = \int_{0}^{t_f} \left[ 
k_1 \left( \frac{\alpha - \bar{\alpha}}{\alpha_{\max}} \right)^2 + 
k_2 \left( \frac{u_{\alpha}}{u_{\alpha,\max}} \right)^2 +
k_3 \left( \frac{u_{\sigma}}{u_{\sigma,\max}} \right)^2
\right] dt,
\end{equation}
where $\alpha$ is the angle of attack, $u_{\alpha}$ is the angle of attack rate, $u_{\sigma}$ is the bank angle rate, $t_f$ is the terminal time, $(k_1, k_2, k_3)$ are design parameters (and are constants), $\bar{\alpha}$ is the value of $\alpha$ where the lift-drag ratio is a maximum, and $(\alpha_{\max},~u_{\alpha,\max},~u_{\sigma,\max})$ are constants. It is noted that $u_{\alpha}$ and $u_{\sigma}$ are the controls for this example.  Next, the dynamics are given as
\begin{equation}\label{eq:Ex-EOM}
\begin{array}{rclcrcl}
\dot{r} 
&=& 
v \left( 1 - 2(\e{B2}^2 + \e{B3}^2) \right),

&\quad&

\dot{v} 
&=& 
\frac{1}{m} \tilde{f}_1, \\

\edot{A1} 
&=&
-\frac{1}{2} \w{A2} \e{A3} + \frac{1}{2} \w{A3} \e{A2},

&\quad&

\edot{B1} 
&=& 
-\frac{1}{2} \w{B2} \e{B3} + \frac{1}{2} \w{B3} \e{B2}, \\

\edot{A2} 
&=& 
\phantom{-} \frac{1}{2} \w{A2} \eta_A - \frac{1}{2} \w{A3} \e{A1},

&\quad&

\edot{B2} 
&=& 
\phantom{-} \frac{1}{2} \w{B2} \eta_B - \frac{1}{2} \w{B3} \e{B1}, \\

\edot{A3} 
&=& 
\phantom{-} \frac{1}{2} \w{A2} \e{A1} + \frac{1}{2} \w{A3} \eta_A,

&\quad&

\edot{B3} 
&=& 
\phantom{-} \frac{1}{2} \w{B2} \e{B1} + \frac{1}{2} \w{B3} \eta_B, \\

\dot{\eta}_A 
&=& 
-\frac{1}{2} \w{A2} \e{A2} - \frac{1}{2} \w{A3} \e{A3},

&\quad&

\dot{\eta}_B 
&=& 
-\frac{1}{2} \w{B2} \e{B2} - \frac{1}{2} \w{B3} \e{B3},
\end{array}
\end{equation}
where $m$ is the mass, $r$ is the geocentric radius, $\{ \e{A1}, \e{A2},\e{A3}, \eta_A \}$ are the $\CAE$ Euler parameters, $v$ is the Earth-relative speed, and $\{ \e{B1}, \e{B2},\e{B3}, \eta_B \}$ are the $\CBA$ Euler parameters.  Note that frame $\C{N}$ is the Earth-centered inertial (ECI) frame, frame $\C{E}$ is the Earth-centered, Earth-fixed (ECEF) frame, and frames $\C{A}$ and $\C{B}$ are the position and velocity frames as defined in Section~\ref{sect:Derivation}.  Moreover, the angular velocity terms in Eq.~\eqref{eq:Ex-EOM} are given by
\begin{equation}\label{eq:Ex-EOM-w}
\begin{array}{rcl}
\w{A2} &=& \frac{2 v}{r} \left( \eta_B \e{B2} - \e{B1} \e{B3} \right), \\

\w{A3} &=& \frac{2 v}{r} \left( \eta_B \e{B3} + \e{B1} \e{B2} \right), \\

\w{B2} &=& - \frac{1}{m v} \tilde{f}_3 
- \w{A2} \left( 1 - 2(\e{B1}^2 + \e{B3}^2) \right)
- 2 \w{A3} (\e{B2} \e{B3} + \e{B1} \eta_B), \\

\w{B3} &=& \frac{1}{m v} \tilde{f}_2 
- 2 \w{A2} (\e{B2} \e{B3} - \e{B1} \eta_B)
- \w{A3} \left( 1 - 2(\e{B1}^2 + \e{B2}^2) \right).
\end{array}
\end{equation}
The apparent force terms in Eqs.~\eqref{eq:Ex-EOM} and \eqref{eq:Ex-EOM-w} are given by Eq.~\eqref{eq:F-apparent-3},
%\begin{equation}\label{eq:Ex-EOM-Fapp}
%\begin{array}{rcl}
%\frac{1}{m} \begin{bmatrix}
%\tilde{f}_1 \\ \tilde{f}_2 \\ \tilde{f}_3
%\end{bmatrix}
%&=&
%\frac{1}{m} \begin{bmatrix}
%f_1 \\ f_2 \\ f_3
%\end{bmatrix}
%
%- 2 \w{e} v
%\begin{bmatrix}
%0 \\ \phantom{-}\CBE(3,3) \\ -\CBE(2,3)
%\end{bmatrix}
%
%- r \w{e}^2 \CBA 
%\begin{bmatrix}
%-\left(\CAE(2,3)\right)^2 - \left(\CAE(3,3)\right)^2 \\
%\CAE(2,3) \CAE(1,3) \\
%\CAE(3,3) \CAE(1,3)
%\end{bmatrix},
%\end{array}
%\end{equation}
where $\w{e}$ now represents the rotation rate of the Earth.  The forces present are lift, drag, and gravity.  Thus, Eq.~\eqref{eq:netForce-B} simplifies to
\begin{equation}\label{eq:Ex-EOM-Ftot}
\begin{array}{rcl}
\begin{bmatrix}
{f}_1 \\ {f}_2 \\ {f}_3
\end{bmatrix}
&=&
\begin{bmatrix}
-D \\ L \cos \sigma \\ L \sin \sigma
\end{bmatrix}

-m \frac{\mu_{e}}{r^2} \CBA(:,1),
\end{array}
\end{equation}
where $L=q S C_L(\alpha)$ is the magnitude of the lift force, $D = q S C_D(\alpha)$ is the magnitude of the drag force, $q=\rho v^2/2$ is the dynamic pressure, $S$ is the reference area, $\sigma$ is the bank angle (defined as a rotation about $\h{b}_1$ from the $\h{b}_2$ direction to the positive lift direction), and $\mu_{e}$ is the gravitational parameter of the Earth.  The details of the aerodynamic model are omitted here but can be found in Refs.~\cite{Clarke2,Clarke1}.

Next, $\alpha$ and $\sigma$ are augmented components of the state such that
\begin{equation}
  \begin{array}{lclclcl}
    \dot{\alpha}& =& u_{\alpha} & \textrm{and} &     \dot{\sigma}& =& u_{\sigma}.
  \end{array}
\end{equation}
Furthermore, the following path constraints are imposed during flight (where $ a_{\max}$, $q_{\min}$, $u_{\alpha,\max}$, and $u_{\sigma,\max}$ are constants):
\begin{equation}\label{eq:Ex-path}
\begin{array}{rclcrcl}
  \sqrt{L^2+D^2}/m & \leq & a_{\max} & , & q & \geq & q_{\min}, \\
  \alpha & \geq & 0 & , & \alpha & \leq & \alpha_{\max}, \\
  |u_{\alpha}| & \leq & u_{\alpha,\max}  & , & |u_{\sigma}| & \leq & u_{\sigma,\max}.
\end{array}
\end{equation}
Next, the boundary conditions are given in Table~\ref{tab:Ex-BCs}.  Note that requiring $\e{B1}(t_f) = \eta_B(t_f) = 0$ imposes the vertical impact condition by forcing the velocity direction ($\h{b}_1$) at impact to point in the opposite direction as the position direction ($\h{a}_1$) at impact.  Finally, the endpoint constraints
\begin{equation}\label{eq:Ex-endpoint}
\begin{array}{lcl}
\textrm{atan2} \left( \CAE(1,2)~~,~~\CAE(1,1) \right) \bigr \rvert_{t_f} & = & \phi_T, \\
\textrm{atan2} \left( \CAE(1,3)~~,~~\sqrt{\left(\CAE(1,1)\right)^2 + \left(\CAE(1,2)\right)^2} \right) \biggr \rvert_{t_f} & = & \theta_T,
\end{array}
\end{equation}
force the terminal position to coincide with the Earth-relative longitude and geocentric latitude of the target, where $\phi_T = 25.15\deg$ is the Earth-relative longitude of the target and $\theta_T = 0\deg$ is the geocentric latitude of the target.

\begin{table}[h]
  \centering
  \caption{Boundary Conditions for the Atmospheric Entry Example.\label{tab:Ex-BCs}}
  \renewcommand{\baselinestretch}{1}\normalsize\normalfont
 \begin{tabular}{|c|c||c|c|}\hline
    Symbol & Value & Symbol & Value \\\hline
    $t_0$ & $0~\textrm{s}$ & $t_f$ & FREE \\
    $r(t_0)$         & $R_{e} + 37~\textrm{km}$ & $r(t_f)$ & $R_{e}$ \\
    $\e{A1}(t_0)$ & $0$ & $\e{A1}(t_f)$ & FREE \\
    $\e{A2}(t_0)$ & $0$ & $\e{A2}(t_f)$ & FREE \\
    $\e{A3}(t_0)$ & $0$ & $\e{A3}(t_f)$ & FREE \\
    $\eta_A(t_0)$ & $1$ & $\eta_A(t_f)$ & FREE \\
    $v(t_0)$         & $7.138~\textrm{km/s}$ & $v(t_f)$ & $1.219~\textrm{km/s}$ \\
    $\e{B1}(t_0)$ & $\sqrt{2} / 2$ & $\e{B1}(t_f)$ & $0$ \\
    $\e{B2}(t_0)$ & $\sqrt{2} / 2$ & $\e{B2}(t_f)$ & FREE \\
    $\e{B3}(t_0)$ & $0$ & $\e{B3}(t_f)$ & FREE \\
    $\eta_B(t_0)$ & $0$ & $\eta_B(t_f)$ & $0$ \\
    $\alpha(t_0)$ & FREE & $\alpha(t_f)$ & $0\deg$ \\
    $\sigma(t_0)$ & FREE & $\sigma(t_f)$ & FREE \\\hline
  \end{tabular}
\end{table}

The atmospheric entry optimal control problem is summarized as follows.  Determine the state $( r(t), $ $\e{A1}(t), \e{A2}(t),$ $\e{A3}(t), \eta_A(t), v(t), \e{B1}(t), \e{B2}(t), \e{B3}(t), \eta_B(t), \alpha(t), \sigma(t) )$ and the control $\left( u_{\alpha}(t), u_{\sigma}(t) \right)$ on the time interval $t \in [0,t_f]$ which minimizes the cost functional of Eq.~\eqref{eq:Ex-Cost} while satisfying the state dynamics, % given by Eqs.~\eqref{eq:Ex-EOM} and \eqref{eq:Ex-EOM-attitude},
the path constraints,
% of Eq.~\eqref{eq:Ex-path},
the endpoint constraints,
% of Eq.~\eqref{eq:Ex-endpoint},
and the boundary conditions.
%listed in Table~\ref{tab:Ex-BCs}.  % Finally, the values of the constants employed in the problem are summarized in Table~\ref{tab:Ex-constants}.

% \begin{table}[h]
%   \centering
%   \caption{Numerical Values for the Constants in the Atmospheric Entry Example.\label{tab:Ex-constants}}
%   \renewcommand{\baselinestretch}{1}\normalsize\normalfont
%  \begin{tabular}{|c|c||c|c|}\hline
%     Symbol & Value & Symbol & Value \\\hline
%     $a_{\max}$ & $0.441~\textrm{km} / \textrm{s}^2$ &
%     $C_{D0}$ & $0.043$ \\
%     $C_{L,\alpha}$ & $1$ &
%     $H$ & $6.914~\textrm{km}$ \\
%     $K$ & $1$ &
%     $k_1$ & $1$ \\
%     $k_2$ & $18$ &
%     $k_3$ & $20$ \\
%     $m$ & $687~\textrm{kg}$ &
%     $q_{\min}$ & $11.97~\textrm{kPa}$ \\
%     $R_{e}$ & $6378.145~\textrm{km}$ &
%     $S$ & $0.6~\textrm{m}^2$ \\
%     $u_{\alpha,\max}$ & $10\deg / \textrm{s}$ &
%     $u_{\sigma,\max}$ & $30\deg / \textrm{s}$ \\

%     $\alpha_{\max}$ & $25\deg$ &
%     $\bar{\alpha}$ & $11.860\deg$ \\
%     $\mu_{e}$ & $3.986 \times 10^5~\textrm{km}^3 / \textrm{s}^2$ &
%     $\rho_0$ & $1.226 \times 10^3~\textrm{kg} / (\textrm{m}^2 \textrm{km})$ \\
%     $\w{e}$ & $7.292 \times 10^{-5}~\textrm{rad} / \textrm{s}$ &
%      &  \\\hline
%   \end{tabular}
% \end{table}

\subsection{Solution Method and Implementation\label{sect:Ex-Implementation}}

The atmospheric entry problem described in Section~\ref{sect:Ex-OCP} is implemented and solved using $hp$-adaptive Gaussian quadrature collocation via the software of Ref.~\cite{Patterson2014}.  Special precautions must be taken when implementing Eq.~\eqref{eq:Ex-EOM} using a collocation method.  Recall that the unit norm constraints $\e{A1}^2 + \e{A2}^2 + \e{A3}^2 + \eta_A^2 = 1$ and $\e{B1}^2 + \e{B2}^2 + \e{B3}^2 + \eta_B^2 = 1$, as well as the angular velocity constraints $\w{A1} = 0$ and $\w{B1} = 0$ are all implicit in Eq.~\eqref{eq:Ex-EOM}.  Thus, the system of differential equations loses four degrees of freedom and may become inconsistent when solving the problem numerically.  The remedy employed here is to introduce four additional control variables, denoted ($u_1$,$u_2$,$u_3$,$u_4$), such that $\edot{A2}  =  u_1$, $ \edot{A3} =  u_2$, $\edot{B2}  =  u_3$, $\edot{B3}  = u_4$ replaces the corresponding differential equations in Eq.~\eqref{eq:Ex-EOM}.  Next, the additional path constraints
\begin{equation}
\begin{array}{rclcrcl}
| u_1 - \left(
\frac{1}{2} \w{A2} \eta_A - \frac{1}{2} \w{A3} \e{A1} \right) | & \leq & \kappa

&,&

| u_2 - \left(
\frac{1}{2} \w{A2} \e{A1} + \frac{1}{2} \w{A3} \eta_A \right) | & \leq & \kappa, \\

| u_3 - \left(
\frac{1}{2} \w{B2} \eta_B - \frac{1}{2} \w{B3} \e{B1} \right) | & \leq & \kappa

&,&

| u_4 - \left(
\frac{1}{2} \w{B2} \e{B1} + \frac{1}{2} \w{B3} \eta_B \right) | & \leq & \kappa,
\end{array}
\end{equation}
are enforced with $\kappa = 10^{-6}$ chosen to be one order of magnitude larger than the NLP solver accuracy tolerance.

\subsubsection{Optimal Entry Trajectory}\label{sect:Ex-solution}

\begin{figure}[hbt!]
  \centering
  \begin{tabular}{lr}
  \subfloat[Altitude, $h(t)$ vs. time, $t$.]{\includegraphics[width=.475\textwidth]{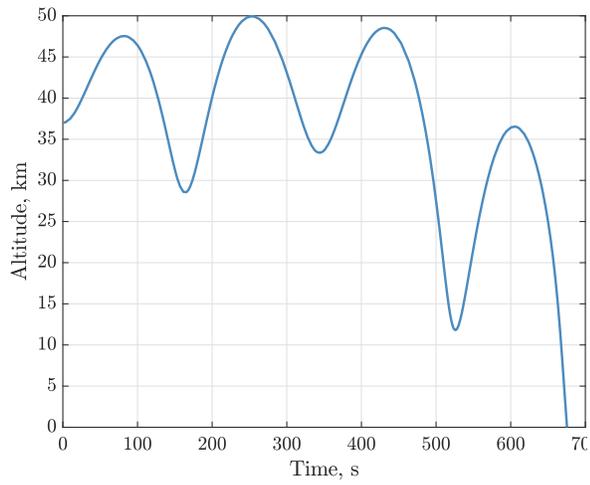}}
  &
  \subfloat[Earth-relative speed, $v(t)$ vs. time, $t$.]{\includegraphics[width=.475\textwidth]{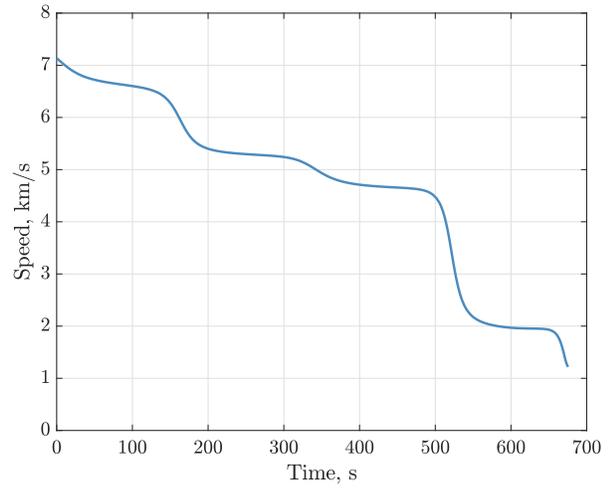}}
  
  \\
  
  \subfloat[$\CAE$ Euler parameters vs. time, $t$.]{\includegraphics[width=.475\textwidth]{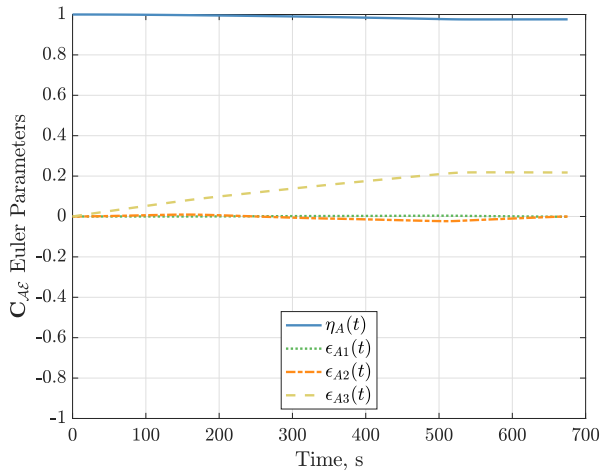}}
  &
  \subfloat[$\CBA$ Euler parameters vs. time, $t$.]{\includegraphics[width=.475\textwidth]{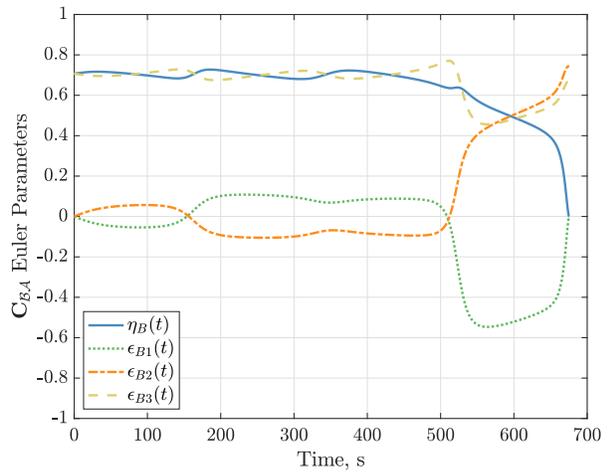}}\label{fig:Ex2-sol-1-CBA}
  
  \\
  
  \subfloat[Angle of attack, $\alpha(t)$ vs. time, $t$.]{\includegraphics[width=.475\textwidth]{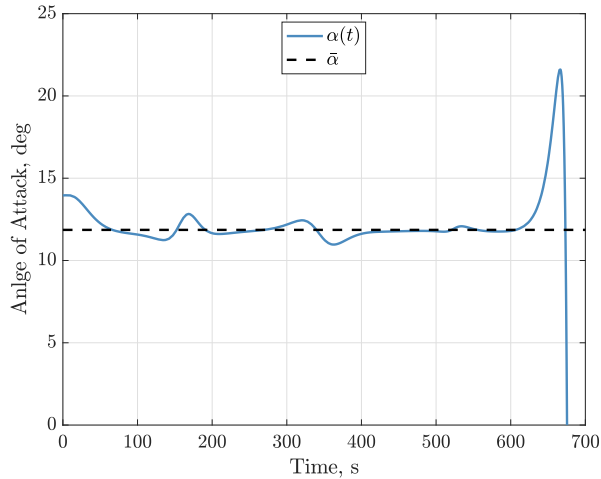}}
  &
  \subfloat[Bank angle, $\{ \sigma(t), \beta(t) \}$ vs. time, $t$]{\includegraphics[width=.475\textwidth]{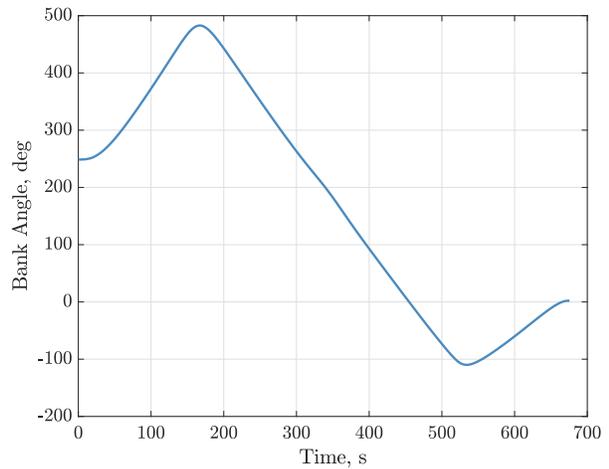}}
  \end{tabular}
  \caption{Optimal solution obtained for the atmospheric entry example.  (f) includes the bank angle $\beta(t)$ of Eq.~\eqref{eq:BS-2} for reference.  Note that $\beta(t_f)$ is undefined and is returned as zero by default.}\label{fig:Ex-sol-1}
\end{figure}

Figure~\ref{fig:Ex-sol-1} shows the optimal entry trajectory where it is seen that the terminal state corresponds to vertically downward flight ($\e{B1}(t_f) = \eta_B(t_f) = 0$).  Notice that Eqs.~\eqref{eq:Ex-EOM} and \eqref{eq:Ex-EOM-w} remain well defined at $\e{B1} = \eta_B = 0$.  In contrast, in a spherical coordinate parameterization the differential equation for azimuth, given by
\begin{equation}\label{eq:Ex-azidot}
\begin{array}{rcl}
\dot{\psi} & =& \frac{L \sin\beta}{m v \cos\gamma}
+ \frac{v}{r} \cos\gamma \sin\psi \tan\theta
- 2 \w{e} \left( \tan\gamma \cos\psi \cos\theta - \sin\theta \right)
+ \frac{r \w{e}^2}{v \cos\gamma} \sin\psi \sin\theta \cos\theta,
\end{array}
\end{equation}
where $\gamma$ is the Earth-relative flight path angle, is not well-defined for vertical flight ($\gamma = \pm 90$ deg).  Thus, the vertical flight singularity using spherical coordinates is eliminated when using the $rv$-Euler parameterization.

%----------------------------------------------
\section{Discussion\label{sect:Discussion}}
The atmospheric entry example in Section~\ref{sect:Example} demonstrates the utility of the parameterizations developed in Section~\ref{sect:Derivation}.  In particular, the example focuses on the ability to use the $rv$-Euler parameters to reformulate the trajectory optimization problems in Refs.~\cite{Clarke2,Clarke1} into a form that does not contain a singularity in vertical flight.  Thus, vertical flight applications such as rocket launches or missile impacts may benefit from using the $rv$- or the $rvL$-Euler parameterizations to model the motion.  Likewise, other applications (such as flight over the poles) may also benefit by modeling the motion with the $rv$-, $rvL$-, or $rvh$-Euler parameterizations due to their lack of singularities.

Each of the three parameterizations developed have distinct advantages over each other as well.  For instance, the $rvh$-Euler parameterization employs just $8$ variables $\{ r, \e{A1}, \e{A2}, \e{A3}, \eta_A, v, \e{B3}, \eta_B \}$ to parameterize position and velocity (as opposed to $10$ in the $rv$- and $rvL$-Euler parameterizations).  The equations of motion are also simpler than those of the $rv$- and $rvL$-Euler parameterizations, but do contain a singularity at $^{\C{E}}\m{h} = \m{0}$ (vertical flight).  The $rv$- and $rvL$-Euler parameterizations have their advantages too, with nonsingular equations of motion in vertical flight at the top of the list.  In addition, the $rvL$-Euler parameterization removes the bank angle from the equations of motion entirely, simplifying the equations slightly and reducing the number of variables needed to model the motion.  Last, the $rv$-Euler parameterization decouples the bank angle rotation from the $\CBA$ Euler parameters, potentially offering slower-varying rates of change for the $\CBA$ Euler parameters when compared to those of the $rvL$-Euler parameterization.

Taken together, the $rv$-, $rvL$-, and $rvh$-Euler parameterizations each offer certain advantages over the others depending upon the particular task at hand.  No doubt parameterizations already in existence offer some advantages over the parameterizations developed in this research as well.  Many considerations must be taken into account when considering employing one parameterization over another.  The presence or absence of singularities, the ease with which forces can be modeled, computational efficiency, scaling, convenience, and physical meaning are a few of the considerations that come to mind.  The relative importance of each consideration will likely vary depending upon the application at hand.  While the $rv$-, $rvL$-, and $rvh$-Euler parameterizations may not be the perfect fit for every application, they do lend themselves as useful alternatives to many parameterizations in wide use today.

%In addition to the absence of singularities, the $rv$-, $rvL$-, and $rvh$-Euler parameterizations greatly reduce the number of trigonometric functions present in the equations of motion when compared to a typical spherical coordinate parameterization.  The significance of avoiding the use of trigonometric functions is related to the computational expense in evaluating a trigonometric function.   For example, consider the following software implementation of $\sin x$ \cite{netlib}, given by
% \begin{equation}\label{eq:sin(x)}
% \sin x \approx x + x z (S_1 + z (S_2 + z (S_3 + z (S_4 + z (S_5 +z S_6))))),
% \end{equation}
%where $z = x^2$ and $S_i,~i = 1,\ldots,6$ are the coefficients of the truncated power series approximation.  The sine evaluation in Eq.~\eqref{eq:sin(x)} requires $14$ floating-point operations, not to mention the function calls, decision tree, etc. which are also executed in the code.  Now, observe in Eqs.~\eqref{eq:Ex1-EOM-Geo} and \eqref{eq:Ex1-EOM-rvEuler} from Section~\ref{sect:Ex1} that $\sin \gamma$ in the spherical parameterization is equivalently expressed as $(1 - 2(\e{B2}^2 + \e{B3}^2))$ in the $rv$-Euler parameterization.  The reduction in the number of floating-point operations is clear.

%----------------------------------------------
\section{Conclusions\label{sect:Conclusion}}
Three parameterizations of the equations of motion for a point mass in flight about a central, rotating body have been derived.  The equations of motion have been shown to avoid singularities found in commonly used parameterizations.  Moreover, it was found that both position-dependent and velocity-dependent forces could be quantified in a tractable manner using any of the three parameterizations.  Finally, a trajectory optimization example involving vertical flight was used to demonstrate the nonsingular nature of the equations of motion.

%----------------------------------------------
\section*{Acknowledgments}
The authors gratefully acknowledge support for this research from the from the U.S.~National Science Foundation under grants CMMI-1563225, DMS-1522629, DMS-1819002, and CMMI‐2031213, from the U.S.~Office of Naval Research under grant N00014-19-1-2543, and from the U.S.~Department of Defense under the National Defense Science \& Engineering Graduate Fellowship (NDSEG) Program.

%----------------------------------------------
\renewcommand{\baselinestretch}{1}
\normalsize\normalfont

\bibliographystyle{aiaa}     % Number the references.
% \bibliography{references}   % Use references.bib to resolve the labels.

\renewcommand{\baselinestretch}{1.5}
\normalsize\normalfont

\end{document}